\documentclass[12pt]{amsart}
\usepackage{euscript,amsmath, amssymb, amsfonts,amsthm}

\theoremstyle{plain}
\newtheorem{thm}{Theorem}[section]
\newtheorem{cor}[thm]{Corollary}
\newtheorem{lem}[thm]{Lemma}
\newtheorem{hypo}[thm]{Hypothesis}

\theoremstyle{definition}
\newtheorem{defn}[thm]{Definition}

\newtheorem{rem}[thm]{Remark}

\makeatletter
\@addtoreset{thm}{section}
\@addtoreset{equation}{section}
\def\theequation{\arabic{section}.\arabic{equation}}
\makeatother

\newcommand{\dv}[2]{\langle #1,#2\rangle}
\newcommand{\Rset}{\mathbb R}
\newcommand{\Cset}{\mathbb C}
\newcommand{\Nset}{\mathbb N}
\renewcommand{\d}{\mathrm{d}}

\begin{document}

\title{Fourier transformation of Sato's hyperfunctions}

\author[A.~G.~Smirnov]{A.~G.~Smirnov$^1$}
\setcounter{footnote}{1}
\footnotetext{The research is partially supported by
the grants RFBR~02-01-00556, INTAS~03-51-6346, and LSS-1578.2003.2.}

\address{I.~E.~Tamm Theory Department, P.~N.~Lebedev
Physical Institute,
Leninsky prospect 53, Moscow 119991, Russia}

\email{smirnov@lpi.ru}

\begin{abstract}
A new generalized function space in which all Gelfand--Shilov classes
$S^{\prime 0}_\alpha$ ($\alpha>1$) of analytic functionals are
embedded is introduced. This space of \emph{ultrafunctionals} does
not possess a natural nontrivial topology and cannot be obtained via
duality from any test function space. A canonical isomorphism between
the spaces of hyperfunctions and ultrafunctionals on $\Rset^k$ is
constructed that extends the Fourier transformation of Roumieu-type
ultradistributions and is naturally interpreted as the Fourier
transformation of hyperfunctions. The notion of carrier cone that
replaces the notion of support of a generalized function for
ultrafunctionals is proposed. A Paley--Wiener--Schwartz-type theorem
describing the Laplace transformation of ultrafunctionals carried by
proper convex closed cones is obtained and the connection between the
Laplace and Fourier transformations is established.
\end{abstract}



\maketitle

\setcounter{footnote}{1}

\section{Introduction}
\label{s1}

It is well known that Sato's hyperfunctions cannot be interpreted as
continuous linear functionals on any test function space. For this
reason, the standard definition of the Fourier transformation of
generalized functions is inapplicable to hyperfunctions. This
difficulty does not appear in the framework of Fourier
hyperfunctions~\cite{Sato} that grow at infinity no faster than any
linear exponential. Kawai~\cite{Kawai} has established that the space
of Fourier hyperfunctions on $\Rset^k$ is naturally identified with
the continuous dual of a suitable test function space (actually
coinciding with the Gelfand--Shilov space $S^1_1(\Rset^k)$) and is
taken to itself by the Fourier transformation. However, the question
is still open whether it is possible to construct the Fourier
transformation of general hyperfunctions with no growth restrictions
imposed. The aim of this paper is to fill this gap.

The proposed construction naturally arises from the consideration of
analytic functionals defined on Gelfand--Shilov test function spaces
$S^0_\alpha(\Rset^k)$ with $\alpha>1$. According to~\cite{GS} the
Fourier transformation induces a topological isomorphism between
$S^0_\alpha(\Rset^k)$ and the space $S^\alpha_0(\Rset^k)$, whose
continuous dual $S^{\prime\alpha}_0(\Rset^k)$ is exactly the space of
Roumieu's ultradistributions~\cite{R} of class $\{k^{\alpha k}\}$.
The space $\mathcal B(\Rset^k)$ of hyperfunctions on $\Rset^k$ can be
thought of as the ``limiting case'' of the spaces
$S^{\prime\alpha}_0(\Rset^k)$ as $\alpha\downarrow 1$. Therefore, we
can try to define the Fourier transform $\mathcal U(\Rset^k)$ of the
space $\mathcal B(\Rset^k)$ by passing to the limit $\alpha\downarrow
1$ in the definition of the spaces $S^{\prime 0}_\alpha(\Rset^k)$.

Unfortunately, we cannot just set $\mathcal U(\Rset^k)=S^{\prime
0}_1(\Rset^k)$ because the space $S^0_1(\Rset^k)$ is
trivial~\cite{GS}. The way of overcoming this difficulty is suggested
by the results of the papers~\cite{Sol1,Sol2} concerning the
localization of analytic functionals belonging to $S^{\prime
0}_\alpha(\Rset^k)$. In these works, the notion of \emph{carrier
cone} that replaces the notion of support of a generalized function
for analytic functionals was proposed (the standard definition of
support does not work because of the lack of test functions with
compact support). The definition of carrier cones is based on
introducing, for every closed cone $K$, a suitable test function
space $S^0_\alpha(K)$ in which $S^{0}_\alpha(\Rset^k)$ is densely
embedded (the precise definition will be given in Section~\ref{s2});
a functional $u\in S^{\prime 0}_\alpha(\Rset^k)$ is said to be
carried by a closed cone $K$ if $u$ has a continuous extension to
$S^0_\alpha(K)$. As shown in~\cite{Sol1}, every functional in
$S^{\prime 0}_\alpha(\Rset^k)$ has a uniquely defined minimal carrier
cone. The definition of the spaces associated with cones is naturally
extended to the case $\alpha=1$ and it turns out that the spaces
$S^0_1(K)$ over \emph{proper}\footnote{A cone $U$ in $\Rset^k$ will
be called proper if $\bar U\setminus\{0\}$ is contained in an open
half-space of $\Rset^k$ (the bar denotes closure). For convex closed
cones, this definition is equivalent to the usual one according to
which a cone is called proper if it contains no straight lines. }
closed cones are nontrivial. The space $\mathcal U(\Rset^k)$ is
obtained by ``gluing together'' the generalized function spaces
$S^{\prime 0}_1(K)$ associated with proper closed cones $K\subset
\Rset^k$ (this procedure will be given a precise meaning in
Section~\ref{s3}).

The properties of the elements of $\mathcal U(\Rset^k)$, which will
be named \emph{ultrafunctionals}, are quite similar to those of
analytic functionals in $S^{\prime 0}_\alpha(\Rset^k)$. In
particular, the definition of carrier cones is extended to the case
of the space $\mathcal U(\Rset^k)$ and it turns out that every
ultrafunctional has a uniquely defined minimal carrier cone. For a
closed proper cone $K$, the space $\mathcal U(K)$ consisting of
ultrafunctionals carried by $K$ coincides with $S^{\prime 0}_1(K)$.
The spaces $S^{\prime 0}_\alpha(K)$ are naturally embedded in
$\mathcal U(K)$ for any closed cone $K$. If $K,K_1,\ldots, K_n$ are
closed cones in $\Rset^k$ such that $K=\bigcup_{j=1}^n K_j$, then
every ultrafunctional $u\in \mathcal U(K)$ is representable in the
form
\begin{equation}\label{0000}
u=\sum_{j=1}^n u_{j},\quad u_j\in \mathcal U(K_j).
\end{equation}
Every exponential decreasing in an open half-space containing a
convex proper closed cone $K$ belongs to the space $S^{0}_1(K)$. This
allows us to define the Laplace transform $\mathcal L_K u$ of every
ultrafunctional $u$ carried by $K$. We prove an elegant
Paley--Wiener--Schwartz-type theorem asserting that the Laplace
transformation $\mathcal L_K$ induces a topological isomorphism
between $\mathcal U(K)$ and the space of all functions analytic in
the tubular domain $\Rset^k+iV$, where $V$ is the interior of the
dual cone of $K$. The Fourier transform $\mathcal F u$ of an
ultrafunctional $u$ carried by a convex proper closed cone $K$ is by
definition the boundary value in $\mathcal B(\Rset^k)$ of the Laplace
transform of $u$. For a general $u\in \mathcal U(\Rset^k)$, we take a
decomposition of the form~(\ref{0000}), where all $K_j$ are convex
proper closed cones, and set $\mathcal F u=\sum_{j=1}^n\mathcal F
u_j$. The hyperfunction $\mathcal F u$ so defined does not depend on
the chosen decomposition. We prove that the operator $\mathcal F$
maps $\mathcal U(\Rset^k)$ isomorphically onto $\mathcal B(\Rset^k)$
and that its restriction to $S^{\prime 0}_\alpha(\Rset^k)$ coincides
with the ordinary Fourier transformation determined via duality by
the Fourier transformation of test functions.

The paper is organized as follows. In Section~\ref{s2}, we give a
brief exposition of the results of the works~\cite{Sol1,Sol2}
concerning the spaces $S^{\prime 0}_\alpha(\Rset^k)$ with
$\alpha>1$ and obtain a useful representation of $S^{\prime
0}_\alpha(\Rset^k)$ in terms of the spaces associated with proper
closed cones. In Section~\ref{s3}, we introduce the spaces
$S^0_1(K)$ and $\mathcal U(K)$ and give the precise formulations of
the main results. In the same section, we prove the compatibility
of the operator $\mathcal F$ with the Fourier transformation of
ultradistributions. Section~\ref{s4} is devoted to a detailed study
of the spaces $S^0_1(K)$ over proper closed cones and to the proof
of the above-mentioned PWS-type theorem. In Section~\ref{s5}, the
results concerning carrier cones (the existence of a unique minimal
carrier cone of an ultrafunctional and the existence of
decompositions of the form~(\ref{0000})) are established. In
Section~\ref{s6}, the bijectivity of the Fourier operator $\mathcal
F$ is proved. In Section~\ref{s7}, we indicate some possible
further developments of these results. The proofs of some algebraic
statements of Section~\ref{s5} are given in Appendices~A and~B.

\textbf{Acknowledgements.}
The author
is grateful to M.~A.~Soloviev for useful discussions and to
the referee for valuable suggestions.

\section{Localization of analytic functionals on Gelfand--Shilov
spaces}
\label{s2}

The space $S_\alpha^\beta(\Rset^k)$ is by definition~\cite{GS}
the union (inductive limit) with respect to $A,B>0$ of the
Banach spaces consisting of smooth functions on
$\Rset^k$ with the finite norm
\begin{equation}
\sup_{x\in \Rset^k,\,\lambda,\,\mu}
\frac{|x^\lambda\partial^\mu f(x)|}{A^{|\lambda|} B^{|\mu|}
|\lambda|^{\alpha |\lambda|} |\mu|^{\beta|\mu|}},
\label{xx0}
\end{equation}
where $\lambda$ and $\mu$ run over all multi-indices and the standard
multi-index notation is used. The spaces $S_\alpha^\beta$ are
nontrivial if $\alpha+\beta>1$ or if $\alpha,\beta>0$ and
$\alpha+\beta=1$. For $\alpha=0$, the spaces $S_\alpha^\beta$ consist
of functions of compact support. If $0\leq \beta<1$, then
$S_\alpha^\beta$ consists of (the restrictions to $\Rset^k$ of)
entire analytic functions and an alternative description of these
spaces in terms of complex variables is possible~\cite{GS}. Namely,
an analytic function $f$ on $\Cset^k$ belongs to the class
$S_\alpha^\beta$ if and only if
\[
|f(z)|\leq C\exp(-|x/A|^{1/\alpha}+|By|^{1/(1-\beta)}),\quad
z=x+iy\in \Cset^k,
\]
for some $A,B>0$ depending on $f$. For definiteness, we assume the
norm $|\cdot|$ on $\Rset^k$ to be uniform, i.e.,  $|x|=\sup_{1\leq
j\leq k}|x_j|$. As shown in~\cite{GS}, the Fourier transformation
isomorphically maps the space $S_\alpha^\beta$ onto $S_\beta^\alpha$.
The Fourier transformation of generalized functions on
$S^\beta_\alpha$ is defined in a standard way, as the dual mapping of
the Fourier transformation of test functions, and maps
$S^{\prime\beta}_\alpha$ onto $S^{\prime\alpha}_\beta$.

In what follows, we confine our discussion to the case $\beta=0$
which is of primary interest to us, but in fact only the condition
$\beta<1$ guaranteeing the analyticity of test functions is necessary
for the constructions described in the rest of this section.
We say that a cone $W$ is a conic neighborhood of a
cone $U$ if $W$ has an open projection\footnote{The projection
$\Pr W$ of a cone $W\subset \Rset^k$ is by definition the canonical
image of $W$ in the sphere $\mathbb
S_{k-1}=(\Rset^k\setminus\{0\})/\Rset_+$; the projection of $W$ is
meant to be open in the topology of this sphere. Note that the
degenerate cone $\{0\}$ is a cone with an open (empty) projection.}
and contains $U$.

\begin{defn}
\label{dxx1}
Let $\alpha>1$ and $U$ be a nonempty cone in $\Rset^k$. The Banach
space $S^{0,B}_{\alpha,A}(U)$ consists of entire analytic functions
on $\Cset^k$ with the finite norm
\[
\|f\|^\alpha_{U,A,B}=\sup_{z=x+iy\in \Cset^k} |f(z)|
\exp(|x/A|^{1/\alpha}-\delta_U(Bx)-|By|),
\]
where $\delta_U(x)=\inf_{x'\in U}|x-x'|$ is the distance from $x$ to
$U$. The space $S^0_\alpha(U)$ is defined by the relation
\[
S^0_\alpha(U)=\bigcup_{A,B>0,\,W\supset U} S^{0,B}_{\alpha,A}(W),
\]
where $W$ runs over all conic neighborhoods of $U$ and the union is
endowed with the inductive limit topology.
\end{defn}

According to the above, for $U=\Rset^k$, this definition is
equivalent to the initial definition of $S^0_\alpha(\Rset^k)$. From
now on and throughout the paper, all cones in question will be
supposed nonempty. As a rule, the word `nonempty'
will be omitted. In the rest of this section, we assume that the
nontriviality condition $\alpha>1$ is satisfied. If $U'\subset U$,
then the space $S^0_\alpha(U)$ is continuously embedded into
$S^0_\alpha(U')$.
If $W\subset\Rset^k$ is a cone with open projection, then
Definition~\ref{dxx1} gives
\begin{equation}\label{OpPr}
S^0_\alpha(W)=\bigcup_{A,B>0} S^{0,B}_{\alpha,A}(W).
\end{equation}
The following statement is an immediate consequence of
Definition~\ref{dxx1}, formula~(\ref{OpPr}), and the
associativity property of inductive limit topologies.

\begin{lem}
\label{l0000}
Let $U$ be a cone in $\Rset^k$. Then
\begin{equation}\label{OpPr1}
S^0_\alpha(U)=\bigcup_{W\supset U} S^{0}_{\alpha}(W),
\end{equation}
where the union is taken over all conic
neighborhoods of $U$
and is endowed with the inductive
limit topology.
\end{lem}

A closed cone $K$ is called a carrier cone of a
functional $u\in S^{\prime 0}_{\alpha}(\Rset^k)$ if $u$ can be
extended continuously to the space $S_{\alpha}^0(K)$. The following
three basic theorems were established in~\cite{Sol1,Sol2}.

\begin{thm}
\label{txx0}
The space $S_{\alpha}^{0}(\Rset^k)$ is dense in
$S_{\alpha}^{0}(U)$ for any cone $U\subset\Rset^k$.
\end{thm}

\begin{thm}
\label{txx1}
If both $K_1$ and $K_2$ are carrier cones of
$u\in S^{\prime 0}_{\alpha}(\Rset^k)$, then so is $K_1\cap K_2$.
\end{thm}

\begin{thm}
\label{txx2}
Let $K_1$ and $K_2$ be closed cones in $\Rset^k$ and $u\in S^{\prime
0}_{\alpha}(\Rset^k)$ be carried by $K_1\cup K_2$. Then there are
$u_{1,2}\in S^{\prime 0}_{\alpha}(\Rset^k)$ carried by $K_{1,2}$ and
such that $u=u_1+u_2$.
\end{thm}

Theorem~\ref{txx0} shows that the space of the functionals carried by
a closed cone $K$ is naturally identified with the space
$S_{\alpha}^{\prime 0}(K)$. It follows from
Theorem~\ref{txx0} and Lemma~\ref{l0000} that a functional
$u\in S_{\alpha}^{\prime 0}(\Rset^k)$
is carried by a closed cone $K$ if and only if
$u$ has a continuous extension to the space $S_{\alpha}^{0}(W)$
for every conic neighborhood $W$ of $K$.
Theorem~\ref{txx1} implies that
the intersection of an arbitrary family
$\{K_\omega\}_{\omega\in\Omega}$ of carrier cones of a functional
$u\in S^{\prime 0}_{\alpha}(\Rset^k)$ is again a carrier cone of $u$.
Indeed, let $W$ be a conic neighborhood of
$K=\bigcap_{\omega\in\Omega} K_\omega$. Then by standard compactness
arguments, there is a finite family
$\omega_1,\ldots,\omega_n\in\Omega$ such that $\tilde
K=\bigcap_{j=1}^n K_{\omega_j}\subset W$. By Theorem~\ref{txx1},
$\tilde K$ is a carrier cone of $u$ and, therefore, $u$ has a
continuous extension to $S^{0}_{\alpha}(W)$.
Hence $K$ is a carrier cone of $u$. In particular, every functional
$u\in S^{\prime 0}_{\alpha}(\Rset^k)$ has a uniquely defined minimal
carrier cone --- the intersection of all carrier cones of $u$.

\begin{rem}
In \cite{Sol1,Sol2}, only open and closed cones were considered.
The space $S^0_\alpha(W)$ associated with an open cone $W$
was defined by formula~(\ref{OpPr}). For a closed cone $K$,
the space $S^0_\alpha(K)$ was defined as the right-hand side
of (\ref{OpPr1}), where the union is taken over all open cones
$W$ such that $K\setminus \{0\}\subset W$.
Definition~\ref{dxx1} covers both these cases. Using cones with open
projection instead of open cones allows treating
the degenerate cone $\{0\}$ on the same footing as
nondegenerate closed cones. Theorem~\ref{txx0} was actually
proved in~\cite{Sol2} only for open and closed $U$. This
implies that Theorem~\ref{txx0} holds for cones with open
projection and Lemma~\ref{l0000} ensures that it is valid
for arbitrary~$U$.
\end{rem}

Let $K, K'$ be closed cones in $\Rset^k$ such that $K'\subset K$. We
denote by $\rho^\alpha_{K',K}$ the natural mapping from $S^{\prime
0}_\alpha(K')$ to $S^{\prime 0}_\alpha(K)$ (if $u\in S^{\prime
0}_\alpha(K')$ then $\rho^\alpha_{K',K}u$ is the restriction of $u$
to $S^{0}_\alpha(K)$). It follows from Theorems~\ref{txx0},
\ref{txx1}, and~\ref{txx2} that
\begin{enumerate}
\item[a)] The mappings $\rho^{\alpha}_{K',K}$ are injective for
any $K'\subset K$.
\item[b)] If $u\in S^{\prime 0}_\alpha(K_1\cup K_2)$, then there are
$u_{1,2}\in S^{\prime 0}_{\alpha}(K_{1,2})$ such that
$u=\rho^\alpha_{K_1,K_1\cup K_2}u_1+\rho^\alpha_{K_2,K_1\cup
K_2}u_2$.
\item[c)] If $u_{1,2}\in S^{\prime 0}_{\alpha}(K_{1,2})$,
$K_{1,2}\subset K$, and $\rho^\alpha_{K_1,K}u_1=
\rho^\alpha_{K_2,K}u_2$, then there is a $u\in S^{\prime
0}_{\alpha}(K_1\cap K_2)$ such that $u_1 =\rho^\alpha_{K_1\cap
K_2,K_1}u$ and $u_2 =\rho^\alpha_{K_1\cap K_2,K_2}u$.
\end{enumerate}

\begin{rem}
\label{r0}
Starting from the spaces $S^{\prime 0}_\alpha(K)$, one can construct
a flabby sheaf $\mathfrak F_\alpha$ on the sphere $\mathbb
S_{k-1}=(\Rset^k\setminus\{0\})/\Rset_+$. For $Q\subset\mathbb
S_{k-1}$, let $C(Q)$ denote the cone in $\Rset^k$ containing the
origin and such that $\Pr C(Q)=Q$. For an open set $O\subset \mathbb
S_{k-1}$, set $\mathfrak F_\alpha(O)=S^{\prime 0}_\alpha(C(\bar
O))/S^{\prime 0}_\alpha(C(\partial O))$, where $\partial O$ is the
boundary of $O$ and the bar stands for closure in $\mathbb S_{k-1}$.
Proceeding as in Section~9.2 of the book~\cite{Hoermander}, where
hyperfunctions are constructed from analytic functionals, and using
the properties a), b), and c) reformulated in terms of closed subsets
of $\mathbb S_{k-1}$, one can define the restriction mappings
$\mathfrak F_\alpha(O_1)\to \mathfrak F_\alpha(O_2)$ for $O_1\supset
O_2$ and prove that $\mathfrak F_\alpha$ is indeed a flabby sheaf.
Note however that $\mathfrak F_\alpha(\mathbb S_{k-1})= S^{\prime
0}_1(\Rset^k)/S^{\prime 0}_1(\{0\})$. Thus, passing from the spaces
$S^{\prime 0}_\alpha(K)$ to the sheaf $\mathfrak F_\alpha$ leads to
the loss of information concerning the functionals carried by the
origin. Moreover, since $S^{\prime 0}_\alpha(\{0\})$ is dense in
every space $S^{\prime 0}_\alpha(K)$, all information about the
topology of these spaces is also lost.
\end{rem}

\begin{lem}
\label{lxx0}
Let $M$ be a closed cone in $\Rset^k$ and $P$ be a set of closed
subcones of M such that $K_1\cap K_2\in P$ for any $K_1,K_2\in P$.
Suppose there is a finite subset $P'$ of $P$ such that
$M=\bigcup_{K\in P'} K$. Then the space $S^{\prime 0}_\alpha(M)$ is
canonically isomorphic as a topological vector space to
$\varinjlim\nolimits_{K\in P} S^{\prime 0}_\alpha(K)$ (the set $P$ is
meant to be naturally ordered by inclusion).
\end{lem}

The inductive limit in Lemma~\ref{lxx0} is taken, in general, over a
partially ordered but not directed set of indices. The definitions of
the inductive system and inductive limit, which are usually
formulated for the case of a directed set of indices, are immediately
extended to this more general case.
Moreover, the usual inductive limit universality property remains
valid in this more general case. Precise formulations concerning
such generalized inductive systems will be given in
the end of this section.

\begin{proof}[Proof of Lemma~$\ref{lxx0}$.]
For $K\in P$, we denote by $\rho^\alpha_{K}$ and $\iota_K$ the
canonical mapping from $S^{\prime 0}_\alpha(K)$ to
$\varinjlim\nolimits_{K\in  P} S^{\prime 0}_\alpha(K)$ and the
canonical embedding of $S^{\prime 0}_\alpha(K)$ into $\oplus_{K\in
P}S^{\prime 0}_\alpha(K)$ respectively. If $K, K'\in  P$ and
$K'\subset K$, then we have $\rho^\alpha_{K',M}=
\rho^\alpha_{K,M}\rho^\alpha_{K',\,K}$ and by the inductive limit
universality property, there is a unique continuous mapping $l\colon
\varinjlim\nolimits_{K\in  P} S^{\prime 0}_\alpha(K)\to S^{\prime
0}_\alpha(M)$ such that $\rho^\alpha_{K,\,M}=l\rho^\alpha_{K}$ for
any $K\in  P$. It follows from property b) that $l$ is surjective
because $M$ can be represented as a union of finitely many cones
belonging to $ P$. We now prove the injectivity of $l$. Let $N$ be
the subspace of $\oplus_{K\in P}S^{\prime 0}_\alpha(K)$ spanned by
all vectors of the form $\iota_{K'} u-\iota_K\rho^\alpha_{K',K}u$,
where $K,K'\in P$, $K'\subset K$, and $u\in S^{\prime 0}_\alpha(K')$.
The space $\varinjlim\nolimits_{K\in  P} S^{\prime 0}_\alpha(K)$ is
by definition the quotient space $\oplus_{K\in P}S^{\prime
0}_\alpha(K)/N$. It suffices to show that for any $K_1,\ldots,K_n\in
P$ and every $u_1\in S^{\prime 0}_\alpha(K_1),\ldots, u_n\in
S^{\prime 0}_\alpha(K_n)$, the relation
$\rho^\alpha_{K_1,\,M}u_1+\ldots+\rho^\alpha_{K_n,\,M}u_n=0$ implies
that $\iota_{K_1}u_1+\ldots+\iota_{K_n}u_n$ belongs to $N$. The proof
is by induction on $n$. If $n=1$ and $\rho^\alpha_{K_1,M}u_1=0$, then
by a) we have $u_1=0$. We now assume $n>1$ and prove the statement
supposing it holds for $n-1$. Let $K=K_1\cup\ldots\cup K_{n-1}$,
$K'=K\cap K_n$, and
$u=\rho^\alpha_{K_1,\,K}u_1+\ldots+\rho^\alpha_{K_{n-1},\,K}u_{n-1}$.
Then we have $\rho^\alpha_{K_n,\,M}u_n=-\rho^\alpha_{K,\,M}u$ and by
property c), there is a $u'\in S^{\prime 0}_\alpha(K')$ such that
$u_n=\rho^\alpha_{K',\,K_n}u'$ and $u=-\rho^\alpha_{K',\,K}u'$. Let
$K'_j=K_j\cap K_n$, $j=1,\ldots,n-1$. Since $P$ is closed under
finite intersections, we have $K'_j\in P$. By property b), there are
$u'_1\in S^{\prime 0}_\alpha(K'_1),\ldots,u'_{n-1}\in S^{\prime
0}_\alpha(K'_{n-1})$ such that
$u'=\rho^\alpha_{K'_1,\,K'}u'_1+\ldots
+\rho^\alpha_{K'_{n-1},\,K'}u'_{n-1}$.
We therefore obtain
\begin{equation}\label{xx1}
u_n = \rho^\alpha_{K'_1,\,K_n}u'_1+\ldots+
\rho^\alpha_{K'_{n-1},\,K_n}u'_{n-1}.
\end{equation}
Set $v_j=u_j+\rho^\alpha_{K'_j,\,K_j}u'_j$, $j=1,\ldots,n-1$, and
$\nu=\iota_{K_1}v_1+\ldots+\iota_{K_{n-1}}v_{n-1}$. By the induction
hypothesis, we have $\nu\in N$ because
\[
\rho^\alpha_{K_1,\,M}v_1+\ldots+\rho^\alpha_{K_{n-1},\,M}v_{n-1}=
\rho^\alpha_{K_1,\,M}u_1+\ldots+\rho^\alpha_{K_n,\,M}u_n=0.
\]
Further, we have
\begin{multline}
\iota_{K_1}u_1+\ldots+\iota_{K_n}u_n=\nu +
[\iota_{K_n}u_n-\iota_{K'_1}u'_1-\ldots-\iota_{K'_{n-1}}u'_{n-1}]+\\+
(\iota_{K'_1}u'_1-\iota_{K_1}\rho^\alpha_{K'_1,\,K_1}u'_1)+\ldots+
(\iota_{K'_{n-1}}u'_{n-1}-
\iota_{K_{n-1}}\rho^\alpha_{K'_{n-1},\,K_{n-1}}u'_{n-1}).
\nonumber
\end{multline}
By definition of the space $N$, the terms in the round brackets
belong to $N$ and in view of~(\ref{xx1}) the term in the square
brackets also belongs to $N$. Therefore, the expression in the
left-hand side belongs to $N$ and the injectivity of $l$ is proved.
It remains to show that $l^{-1}$ is continuous. Suppose at first that
the set $P$ is finite. Since $S^{\prime 0}_\alpha(K)$ are Fr\'echet
spaces~\cite{Sol2}, $\oplus_{K\in P}S^{\prime 0}_\alpha(K)$ is also a
Fr\'echet space. By the above, $N$ coincides with the kernel of the
continuous mapping $\{u_K\}_{K\in P}\to \sum_{K\in
P}\rho^\alpha_{K,\,M} u_K$. Therefore, $N$ is a closed subspace of
$\oplus_{K\in P}S^{\prime 0}_\alpha(K)$ and
$\varinjlim\nolimits_{K\in  P} S^{\prime 0}_\alpha(K)$ is a Fr\'echet
space. The continuity of $l^{-1}$ now follows from the open mapping
theorem. If $P$ is arbitrary, then choose a finite set $P'$ such that
$M=\bigcup_{K\in P'}K$. We can assume that $P'$ is closed under
intersections of its elements (otherwise we can add to $P'$ all cones
that are intersections of elements of $P'$). Let $l'$ and $m$ be the
canonical mappings from $\varinjlim\nolimits_{K\in  P'} S^{\prime
0}_\alpha(K)$ to $S^{\prime 0}_\alpha(M)$ and from
$\varinjlim\nolimits_{K\in  P'} S^{\prime 0}_\alpha(K)$ to
$\varinjlim\nolimits_{K\in  P} S^{\prime 0}_\alpha(K)$ respectively.
Then we have $l'=lm$. By the above, $l'$ is a topological isomorphism
and, therefore, $l^{-1}=m l'^{-1}$ is continuous. The lemma is
proved.
\end{proof}

In particular, the conditions of Lemma~\ref{lxx0} are satisfied if
$P$ is equal to the set $\mathcal P(M)$ of all nonempty closed proper
subcones of $M$. We thus have the canonical isomorphism
\begin{equation}\label{xx2}
S^{\prime 0}_\alpha(M)\backsimeq
\varinjlim\nolimits_{K\in \mathcal P(M)} S^{\prime 0}_\alpha(K).
\end{equation}

We end this section by reformulating some standard definitions
and facts
related to inductive limits for the case of partially ordered, but
not necessarily directed sets of indices. By an inductive system
$\mathcal X$ of (locally convex topological) vector spaces indexed by
a partially ordered set $A$, we mean the following data:\footnote{ In
the rest of this section and in Section~\ref{s5},
where abstract
inductive systems are discussed, the Greek letter $\alpha$ is
systematically used to denote an element of a partially ordered set
$A$ and has nothing in common with the index of Gelfand--Shilov
spaces.}
\begin{enumerate}
\item a family $\{\mathcal X(\alpha)\}_{\alpha\in A}$ of (locally
convex topological) vector spaces;
\item a family of (continuous) linear mappings $\rho^{\mathcal
X}_{\alpha\alpha'}\colon\mathcal X(\alpha)\to \mathcal X(\alpha')$
defined for $\alpha\leq \alpha'$ and satisfying the conditions
\begin{enumerate}
\item[(i)] $\rho^{\mathcal X}_{\alpha\alpha}$ is the identity mapping
for any $\alpha\in A$;
\item[(ii)] $\rho^{\mathcal X}_{\alpha\alpha''}=\rho^{\mathcal
X}_{\alpha'\alpha''}\rho^{\mathcal X}_{\alpha\alpha'}$ for
$\alpha\leq\alpha'\leq\alpha''$.
\end{enumerate}
\end{enumerate}
In other words, $\mathcal X$ is a covariant functor from the small
category $A$ to the category of (locally convex topological) vector
spaces. Let $\iota^{\mathcal X}_\alpha$ denote the canonical
embedding of $\mathcal X(\alpha)$ in $\oplus_{\alpha'\in A} \mathcal
X(\alpha')$. The inductive limit $\varinjlim_{\alpha\in A}\mathcal
X(\alpha)$ (or simply $\varinjlim\mathcal X$) is by definition the
quotient space $[\oplus_{\alpha\in A} \mathcal X(\alpha)]/N^{\mathcal
X}$, where $N^{\mathcal X}$ is the subspace of $\oplus_{\alpha\in A}
\mathcal X(\alpha)$ spanned by all elements of the form
$\iota^{\mathcal X}_\alpha x- \iota^{\mathcal
X}_{\alpha'}\rho^{\mathcal X}_{\alpha\alpha'}x$, $x\in \mathcal
X(\alpha)$. The canonical mapping $\rho^{\mathcal X}_\alpha\colon
\mathcal X(\alpha)\to \varinjlim\mathcal X$ is defined by the
relation $\rho^{\mathcal X}_\alpha=j^\mathcal X \iota^{\mathcal
X}_\alpha$, where $j^\mathcal X$ is the canonical surjection of
$\oplus_{\alpha\in A} \mathcal X(\alpha)$ onto $\varinjlim\mathcal
X$. As usual, we have the following inductive limit universality
property:
\begin{enumerate}
\item[] Let $E$ be a (locally convex topological) vector space and
$h_\alpha$ be (continuous) linear mappings from $\mathcal X(\alpha)$
to $E$ such that $h_{\alpha'}\rho^{\mathcal
X}_{\alpha\alpha'}=h_\alpha$ for any $\alpha\leq \alpha'$. Then there
is a unique (continuous) linear mapping $h\colon\varinjlim\mathcal
X\to E$ such that $h_\alpha=h\rho^{\mathcal X}_\alpha$ for any
$\alpha\in A$.
\end{enumerate}

\begin{rem}
\label{r1}
Although the above definitions are quite standard, the properties of
such generalized inductive systems may be very different from those
of inductive systems indexed by directed sets. For example, canonical
mappings $\rho^{\mathcal X}_\alpha$ may be not injective even if all
connecting morphisms $\rho^{\mathcal X}_{\alpha\alpha'}$ are
injective. Indeed, let $E$ be a vector space and $A$ be the
four-element set $\{\alpha,\beta,\gamma,\delta\}$ with the order
defined by the relations $\alpha\leq \gamma$, $\alpha\leq\delta$,
$\beta\leq \gamma$, and $\beta\leq\delta$. We define the inductive
system $\mathcal X$ setting $\mathcal X(\alpha)=\mathcal
X(\beta)=\mathcal X(\gamma)=\mathcal X(\delta)= E$ and
$\rho^{\mathcal X}_{\alpha\gamma}=-\rho^{\mathcal X}_{\alpha\delta}=
\rho^{\mathcal X}_{\beta\gamma}=\rho^{\mathcal
X}_{\beta\delta}=\mathrm{id}_E$, where $\mathrm{id}_E$ is the
identity mapping. Fix $x\in E$ and set $z_1=\iota^{\mathcal
X}_\alpha\,x -\iota^{\mathcal X}_\gamma\,x$, $z_2=\iota^{\mathcal
X}_\alpha\,x +\iota^{\mathcal X}_\delta\,x$, $z_3=\iota^{\mathcal
X}_\beta\,x -\iota^{\mathcal X}_\delta\,x$, and $z_4=\iota^{\mathcal
X}_\gamma\,x -\iota^{\mathcal X}_\beta\,x$. Obviously
$z_1,\ldots,z_4\in N^{\mathcal X}$ and, therefore, $\iota^{\mathcal
X}_\alpha\,x=(z_1+z_2+z_3+z_4)/2$ belongs to $N^{\mathcal X}$. This
means that $\rho^{\mathcal X}_\alpha\,x=0$ for any $x\in E$.
\end{rem}

\section{Basic definitions and formulations of main results}
\label{s3}

We now extend the constructions of the preceding section to the case
$\alpha=1$ which is of primary interest to us. By analogy with
Definition~\ref{dxx1}, we introduce suitable test function spaces
associated with cones in $\Rset^k$.

\begin{defn}
\label{dxx2} Let $U$
be a cone in $\Rset^k$. The Banach space $S^{0,B}_{1,A}(U)$ consists
of entire analytic functions on $\Cset^k$ with the finite norm
\[
\|f\|_{U,A,B}=\sup_{z=x+iy\in \Cset^k}
|f(z)|\exp(|x/A|-\delta_U(Bx)-|By|),
\]
where $\delta_U(x)=\inf_{x'\in U}|x-x'|$ is the distance from $x$ to
$U$.
The space $S^0_1(U)$ is defined by the relation
\[
S^0_1(U)=\bigcup_{A,B>0,\,W\supset U} S^{0,B}_{1,A}(W),
\]
where $W$ runs over all conic neighborhoods of $U$
and the union is endowed with the
inductive limit topology.
\end{defn}

For $U=\Rset^k$, Definition~\ref{dxx2} is equivalent to
the standard definition of $S^0_1(\Rset^k)$ given
in~\cite{GS}. Therefore, the space $S^0_1(U)$ is
trivial for $U=\Rset^k$. A sufficient condition for the
nontriviality of $S^0_1(U)$ will be given in Lemma~\ref{lxx2}.
Representation~(\ref{OpPr}) for the spaces $S^0_\alpha(W)$
associated with cones
with open projection and Lemma~\ref{l0000} obviously remain
valid for $\alpha=1$ (in fact, one can show that
formula~(\ref{OpPr}) holds for $\alpha=1$ even without
the assumption that $W$ has an open projection, but
we shall not prove this fact here).
As shown in~\cite{Sol2}, $S^0_\alpha(U)$ with $\alpha>1$ are
DFS-spaces (we recall that DFS-spaces are, by definition, the
inductive limits of injective compact sequences of locally convex
spaces). In particular, they (and their duals) are reflexive,
complete, and Montel spaces~\cite{Komatsu}. The following lemma shows
that the spaces $S^0_1(U)$ enjoy the same nice topological
properties.

\begin{lem}
\label{lxxx0}
The space $S^0_1(U)$ is DFS for any cone $U\subset\Rset^k$.
\end{lem}

\begin{proof}
It suffices to show that the inclusion mapping
$S^{0,B}_{1,A}(U)\to S^{0,B'}_{1,A'}(U)$ is compact for any $A'>A$,
$B'>B$ and every cone $U\subset\Rset^k$.
Let $\{f_m\}_{m\in\Nset}$ be a sequence of functions
belonging to the unit ball of the space $S^{0,B}_{1,A}(U)$. By
Montel's theorem, this sequence contains a subsequence $\{f_{m_n}\}$
which converges uniformly on compact sets in $\Cset^k$ to an entire
analytic function $f$. To prove the statement, it suffices to show
that the sequence $\{f_{m_n}\}$ converges to $f$ in
$S^{0,B'}_{1,A'}(U)$.
Set $\varrho_{U,A,B}(x+iy)=-|x/A|+B\delta_U(x)+B|y|$,
$x,y\in\Rset^k$.
Since $\|f_{m_n}\|_{U,A,B}\leq 1$, we have
$|f_{m_n}(z)|\leq e^{\varrho_{U,A,B}(z)}$.
Passing to the limit $n\to \infty$ in this
inequality,
we conclude that
$f\in S^{0,B}_{1,A}(U)$ and $\|f\|_{U,A,B}\leq 1$. Further, for any
$R>0$, we have
\begin{multline}
\|f-f_{m_n}\|_{U,A',B'}\leq e^{R/A'}\sup_{|z|\leq R}
|f(z)-f_{m_n}(z)|+\\+
\|f-f_{m_n}\|_{U,A,B}\sup_{|z|>R}e^{\varrho_{U,A,B}(z)-
\varrho_{U,A',B'}(z)}\leq e^{R/A'}\sup_{|z|\leq R}|f(z)-f_{m_n}(z)| +
2e^{-LR},\nonumber
\end{multline}
where $L=\min((A'-A)/A'A, B'-B)$. Choose $R(\varepsilon)$ and
$n(\varepsilon)$ such that $2e^{-LR(\varepsilon)}<\varepsilon/2$ and
$e^{R(\varepsilon)/A'}
\sup_{|z|\leq R(\varepsilon)}|f(z)-f_{m_n}(z)|<\varepsilon/2$ for
any $n\geq n(\varepsilon)$. Then
$\|f-f_{m_n}\|_{U,A',B'}<\varepsilon$ for any $n\geq n(\varepsilon)$.
The lemma is proved.
\end{proof}

\begin{lem}
\label{lxx2}
Let $U$ be a cone in $\Rset^k$. If $U$ is a proper cone, then
the space $S^0_1(U)$ is nontrivial. If $U$ contains a straight line,
then $S^0_1(U)$ is trivial.
\end{lem}

\begin{proof}
Let $U$ be a proper cone and $l$ be a linear functional on $\Rset^k$
such that $\bar U\setminus\{0\}\subset \{\,x\in\Rset^k\,|\,l(x)>0\}$.
Then
$S^0_1(U)$ contains the function $f(z)=e^{-l(z)}$ and, therefore, is
nontrivial. Now let $U$ contain a straight line and $f\in S^0_1(U)$.
Let $W$ be a conic neighborhood of $U$ such that $f\in
S^{0,B}_{1,A}(W)$ for some $A,B>0$ and $\tilde W$ be the union of
all straight lines contained in
$W$. Clearly, $\tilde W$ is a cone with a nonempty interior.
For $x\in \tilde W\setminus\{0\}$
and $\tau\in \Cset$, we set $g(\tau)=f(\tau x)$.
It easily follows from Definition~\ref{dxx2} that
$g\in
S^{0,B|x|}_{1,A/|x|}(\Rset)$ and hence $g\equiv 0$. Therefore,
$f(x)=g(1)=0$, i.e.,
$f$ vanishes on $\tilde W$. By the uniqueness theorem, we
conclude that $f$ is identically zero. The lemma is proved.
\end{proof}

Lemma~\ref{lxx2} suggests that we can try to define the desired
``nontrivialization'' $\mathcal U(\Rset^k)$ of the space $S^{\prime
0}_1(\Rset^k)$ (and, more generally, of the space $S^{\prime 0}_1(M)$
over an arbitrary closed cone $M$) as the right-hand side
of~(\ref{xx2}) with $\alpha=1$. We then arrive at the following
definition.

\begin{defn}
\label{dxx3}
Let $M$ be a closed cone in $\Rset^k$. The space $\mathcal U(M)$ is
defined to be the inductive limit $\varinjlim\nolimits_{K\in \mathcal
P(M)} S^{\prime 0}_1(K)$, where $\mathcal P(M)$ is the set of all
nonempty proper closed cones contained in $M$. The elements of
$\mathcal U(\Rset^k)$ are called ultrafunctionals. A closed cone $K$
is said to be a carrier cone of an ultrafunctional $u$ if the latter
belongs to the image of the canonical mapping from $\mathcal U(K)$ to
$\mathcal U(\Rset^k)$.
\end{defn}

In this definition, the set $\mathcal P(M)$ is meant to be ordered by
inclusion and the inductive limit is taken with respect to the
natural morphisms $\rho_{K',\,K}\colon S^{\prime 0}_1(K')\to
S^{\prime 0}_1(K)$ which are defined for $K'\subset K$ and map the
functionals belonging to $S^{\prime 0}_1(K')$ into their restrictions
to the space $S^{0}_1(K)$. The canonical mappings from $\mathcal
U(K')$ to $\mathcal U(K)$ will be denoted by $\rho^{\mathcal
U}_{K',\,K}$. Note that if $K$ is a proper closed cone, then
${\mathcal U}(K)$ is canonically isomorphic to~$S^{\prime 0}_1(K)$.

We shall see that $\mathcal U(\Rset^k)$ is Fourier-isomorphic to the
space $\mathcal B(\Rset^k)$ which is known to have no natural
topology. Therefore, the following result is by no means surprising.

\begin{lem}
\label{lxx3}
Let $M$ be a closed cone in $\Rset^k$ containing a straight line.
Then the inductive limit topology on ${\mathcal U}(M)$ is trivial
(i.e., ${\mathcal U}(M)$ and $\varnothing$ are the only open sets).
\end{lem}

\begin{proof}
It suffices to prove that any continuous linear functional $l$ on
${\mathcal U}(M)$ is equal to zero. For $K\in \mathcal P(M)$, we
denote by $\rho_K$ the canonical mapping from $S^{\prime 0}_1(K)$ to
${\mathcal U}(M)$. The continuity of $l$ means that the functional
$l\rho_K$ is continuous on $S^{\prime 0}_1(K)$ for any $K\in \mathcal
P(M)$. By Lemma~\ref{lxxx0}, the space $S^0_1(K)$ is reflexive for
any cone $K$. Hence for any $K\in \mathcal P(M)$, there is a function
$f_K\in S^0_1(K)$ such that $l\rho_K\,u=u(f_K)$ for every $u\in
S^{\prime 0}_1(K)$. If $K'\subset K$, then we have
\[
u(f_{K'})=l\rho_{K'}\,u=
l\rho_{K}\rho_{K',\,K}\,u=
(\rho_{K',\,K}u)(f_{K})=u(f_{K}),\quad u\in S^{\prime 0}_1(K'),
\]
and, consequently, $f_{K'}=f_{K}$. Choosing $K'$ equal to the
degenerate cone $\{0\}$, we see that $f_K=f_{\{0\}}$ does not depend
on $K\in\mathcal P(M)$ and, therefore, belongs to the space
$L=\bigcap_{K\in \mathcal P(M)} S^0_1(K)$. Let
$K_1,\ldots,K_n\in\mathcal P(M)$ be such that $M=K_1\cup\ldots\cup
K_n$. Since $\delta_{K_1\cup\ldots\cup
K_n}(x)=\min(\delta_{K_1}(x),\ldots,\delta_{K_n}(x))$ for any $x\in
\Rset^k$, it follows from Definition~\ref{dxx2} that
$S^0_1(M)=S^0_1(K_1\cup\ldots\cup K_n)= S^0_1(K_1)\cap\ldots\cap
S^0_1(K_n)$. Hence $L\subset S^0_1(M)$ is trivial by Lemma~\ref{lxx2}
and $f_K=0$ for any $K\in\mathcal P(M)$. This means that $l\rho_K=0$
for every $K\in\mathcal P(M)$ and, therefore, $l=0$. The lemma is
proved.
\end{proof}

Thus, there is, in general, no natural way to define a reasonable
topology on ${\mathcal U}(K)$. Because of this, we do not endow these
spaces with any topology and consider them only from algebraic point
of view. One of the main results of this paper is that the
ultrafunctionals have the same localization properties as the
analytic functionals belonging to $S^{\prime 0}_{\alpha}(\Rset^k)$
with $\alpha>1$. More precisely, the following analogues of
Theorems~\ref{txx0},~\ref{txx1}, and \ref{txx2} are valid.

\begin{thm}
\label{txx3}
The natural mapping $\rho^{\mathcal U}_{K',\,K}\colon {\mathcal
U}(K')\to{\mathcal U}(K)$ is injective for any closed cones $K$ and
$K'$ such that $K'\subset K$.
\end{thm}

\begin{thm}
\label{txx4}
Let $\{K_\omega\}_{\omega\in\Omega}$ be an arbitrary family of
carrier cones of an ultrafunctional $u$. Then
$\bigcap_{\omega\in\Omega} K_{\omega}$ is also a carrier cone of $u$.
\end{thm}

\begin{thm}
\label{txx5}
Let $K_1$ and $K_2$ be closed cones in $\Rset^k$ and an
ultrafunctional $u$ be carried by $K_1\cup K_2$. Then there are
$u_{1,2}\in {\mathcal U}(\Rset^k)$ carried by $K_{1,2}$ such that
$u=u_1+u_2$.
\end{thm}

These theorems will be proved in Section~\ref{s5}.

\begin{rem}
\label{r0a}
The spaces ${\mathcal U}(K)$ determine a flabby sheaf on the sphere
$\mathbb S_{k-1}$ in the same way as the spaces $S^{\prime
0}_\alpha(K)$ (see Remark~\ref{r0}).
\end{rem}

For $u\in S^{\prime 0}_1(K)$, one can in a standard way define the
operators of partial differentiation and multiplication by an entire
function $g$ of infra-exponential type (i.e., satisfying the bound
$|g(z)|\leq C_\varepsilon e^{\varepsilon|z|}$ for every
$\varepsilon>0$):
\[
\partial u/\partial x_j(f)=-u(\partial f/\partial x_j),\quad (gu)(f)=
u(gf),\quad f\in S^0_1(K),\,j=1,\ldots,k.
\]
These operations are obviously compatible with the connecting
morphisms $\rho_{K',\,K}$ and, therefore, can be lifted to the spaces
${\mathcal U}(K)$ over arbitrary closed cones. Let $\alpha>1$. The
natural mappings from $S^{\prime 0}_\alpha(K)$ to $S^{\prime 0}_1(K)$
taking functionals in $S^{\prime 0}_\alpha(K)$ to their restrictions
to $S^0_1(K)$ are compatible with the connecting morphisms
$\rho^\alpha_{K',\,K}$ and $\rho_{K',\,K}$ and in view of (\ref{xx2})
determine a mapping from $S^{\prime 0}_\alpha(K)$ to ${\mathcal
U}(K)$ for any closed cone $K$. Below we shall see that these
mappings are injective, i.e., the space $S^{\prime 0}_\alpha(K)$ can
be regarded as a subspace of ${\mathcal U}(K)$.

We now describe the construction of the Fourier transformation of
hyperfunctions. As a first step, we consider the Laplace
transformation of analytic functionals on the spaces $S^{0}_1(K)$
over convex proper closed cones. In the rest of this section, we
identify $S^{\prime 0}_1(K)$ with ${\mathcal U}(K)$ for $K\in
\mathcal P(\Rset^k)$. For brevity, the natural embeddings
$\rho^{\mathcal U}_{K,\,\Rset^k}\colon{\mathcal U}(K)\to {\mathcal
U}(\Rset^k)$ and $\rho^\alpha_{K,\,\Rset^k}\colon S^{\prime
0}_\alpha(K)\to S^{\prime 0}_\alpha(\Rset^k)$ will be denoted by
$\sigma_K$ and $\sigma^\alpha_K$ respectively. Let
$\dv{\cdot}{\cdot}$ be a symmetric nondegenerate bilinear form on
$\Rset^k$. Given a cone $U\subset \Rset^k$, we denote by $U^*$ its
dual cone $\{\,x\in \Rset^k\,|\,\dv{x}{\eta}\geq 0 \mbox{ for any }
\eta\in U\,\}$. Note that $U^*$ is always closed and convex. A cone
$U$ is proper if and only if $U^*$ has a nonempty interior. If $V$ is
an open cone, then the function $e^{i\dv{\cdot}{\zeta}}$ belongs to
$S^0_1(V^*)$ for every $\zeta\in
T^V\stackrel{\mathrm{def}}{=}\Rset^k+iV$. Given an open set $O\subset
\Rset^k$, we denote by $\mathcal{A}(O)$ the space of functions
analytic in an open set $T^O\subset\Cset^k$. The space
$\mathcal{A}(O)$ is endowed with the topology of uniform convergence
on compact subsets of $T^O$.

\begin{thm}
\label{txx6}
Let $K$ be a convex proper closed cone in $\Rset^k$ and
$V=\mathrm{int}\,K^*$. For any $u\in {\mathcal U}(K)$, the function
$\zeta\to u(e^{i\dv{\cdot}{\zeta}})$ is analytic in $T^V$. The linear
mapping $\mathcal L_K\colon {\mathcal U}(K)\to \mathcal{A}(V)$ taking
$u\in {\mathcal U}(K)$ to this function is a topological isomorphism.
\end{thm}

The proof of this theorem will be given in Section~\ref{s4}. The
function $\mathcal L_V u$ is called the Laplace transform of $u$. By
definition, we have
\begin{equation}\label{7}
(\mathcal L_K u)(\zeta)=u(e^{i\dv{\cdot}{\zeta}}),\quad \zeta\in
\Rset^k+i\,\mathrm{int}\,K^*.
\end{equation}
For an open cone $V\subset \Rset^k$, we denote by $\mathbf b_V$ the
linear mapping taking functions in $\mathcal{A}(V)$ to their boundary
values in the space of hyperfunctions $\mathcal B(\Rset^k)$. Let
$K,K'\subset \Rset^k$ be proper convex closed cones,
$V=\mathrm{int}\,K^{*}$, and $V'=\mathrm{int}\,K^{\prime *}$. If
$K'\subset K$, then $\mathcal L_{K}\,\rho_{K',\,K}\,u$ is the
restriction of $\mathcal L_{K'}\,u$ to $T^{V}$ for any $u\in
S^{\prime 0}_1(K')$. This implies that $\mathbf b_{V'}L_{K'}= \mathbf
b_{V}\mathcal L_{K}\,\rho^{\mathcal U}_{K',\,K}$, and by the
inductive limit universality property,\footnote{Note that in the
definition of ${\mathcal U}(\Rset^k)$, it suffices to take the
inductive limit over all proper \emph{convex} closed cones in
$\Rset^k$ because the convex hull of a proper closed cone is again a
proper closed cone.} there is a unique mapping $\mathcal
F\colon{\mathcal U}(\Rset^k)\to \mathcal B(\Rset^k)$ such that
\begin{equation}\label{7a}
\mathcal F\sigma_{V^*}=\mathbf b_{V}\mathcal L_{V^*}
\end{equation}
for any open convex cone $V\subset \Rset^k$.

\begin{thm}
\label{txx7}
The operator $\mathcal F$ maps ${\mathcal U}(\Rset^k)$ isomorphically
onto $\mathcal B(\Rset^k)$.
\end{thm}

This theorem will be proved in Section~\ref{s6}. The operator
$\mathcal F$ is naturally interpreted as the inverse Fourier
transformation of hyperfunctions. Indeed, for any $j=1,\ldots,k$ and
$u\in {\mathcal U}(\Rset^k)$, we obviously have the standard
relations
\[
\mathcal F\left[\frac{\partial u}{\partial x_j}\right](\xi)=
-i \xi_j[\mathcal F u](\xi),\quad
\mathcal F\left[x_j u\right](\xi) =
-i\frac{\partial [\mathcal Fu]}{\partial \xi_j}.
\]
Moreover, the restriction of $\mathcal F^{-1}$ to ultradistributions
of the class $S^{\prime \alpha}_0(\Rset^k)$ coincides with the
ordinary Fourier transformation determined via duality by the Fourier
transformation of test functions. More precisely, let $\alpha>1$ and
the Fourier transformation $\hat f$ of a test function $f\in
S^{\alpha}_0(\Rset^k)$ be defined by the relation $\hat f(x)=\int
f(\xi)e^{\dv{x}{\xi}}\,\d\xi$. As mentioned in Section~\ref{s2}, the
mapping $f\to\hat f$ is a topological isomorphism from
$S^{\alpha}_0(\Rset^k)$ onto $S_{\alpha}^0(\Rset^k)$. Let $\mathcal
F^\alpha$ denote its dual mapping acting on generalized functions.
Then we have
\begin{equation}\label{8}
\mathcal F e^\alpha=i^\alpha \mathcal F^\alpha,
\end{equation}
where $e^\alpha\colon S^{\prime 0}_\alpha(\Rset^k)\to {\mathcal
U}(\Rset^k)$ and $i^\alpha\colon S^{\prime \alpha}_0(\Rset^k)\to
\mathcal B(\Rset^k)$ are canonical mappings (see~\cite{Komatsu1} for
the construction of the natural embedding of ultradistributions into
the space of hyperfunctions). To prove~(\ref{8}), we recall some
results concerning the Laplace transformation of analytic functionals
belonging to the spaces $S^{\prime 0}_\alpha$ with $\alpha>1$. For an
open cone $V$, we denote by $\mathcal A^\alpha(V)$ the Fr\'echet
space consisting of functions analytic in $T^V$ and having the finite
norms
\[
|||\mathbf v|||_{V',\varepsilon,R}=
\sup_{\zeta\in T^{V'},\,|\zeta|\leq R}
|\mathbf v(\zeta)|\exp[-\varepsilon|\eta|^{-1/(\alpha-1)}],
\quad \eta=\mathrm{Im}\,\zeta,
\]
for any $\varepsilon, R>0$ and every compact subcone $V'$ of $V$. The
following result has been established in~\cite{Sol2}.

\begin{thm}
\label{txx8}
Let $\alpha>1$, $K$ be a convex proper closed cone in $\Rset^k$, and
$V=\mathrm{int}\,K^*$. For any $u\in S^{\prime 0}_\alpha(K)$, the
function $T^V\ni\zeta\to u(e^{i\dv{\cdot}{\zeta}})$ belongs to
$\mathcal A^\alpha(V)$. The linear mapping $\mathcal L^\alpha_K\colon
S^{\prime 0}_\alpha(K)\to \mathcal A^\alpha(V)$ taking $u\in
S^{\prime 0}_\alpha(K)$ to this function is a topological
isomorphism. The function $(\mathcal L^\alpha_K u)(\cdot+i\eta)$
tends to $\mathcal F^\alpha \sigma^\alpha_K\,u$ in the topology of
$S^{\prime\alpha}_0(\Rset^k)$ as $\eta\to 0$ inside a fixed compact
subcone $V'$ of $V$.
\end{thm}

This theorem implies the existence, for every open convex cone $V$,
of the continuous boundary value operator $\mathbf b^\alpha_V
\colon\mathcal A^\alpha(V)\to S^{\prime\alpha}_0(\Rset^k)$ satisfying
the relation
\begin{equation}\label{9}
\mathcal F^\alpha \sigma^\alpha_{V^*}=
\mathbf b^\alpha_V \mathcal L^{\alpha}_{V^*}.
\end{equation}
Let $j^\alpha_V$ be the inclusion of $\mathcal{A}^\alpha(V)$ into
$\mathcal{A}(V)$ and $e^\alpha_K$ be the canonical mapping from
$S^{\prime 0}_\alpha(K)$ to ${\mathcal U}(K)$ (in particular,
$e^\alpha_{\Rset^k}=e^\alpha$). By definition of the mappings
$e^\alpha_K$, $\mathcal L^\alpha_K$, and $\mathcal L_K$, we have the
relations $j^\alpha_V\mathcal L^\alpha_{V^*}=\mathcal
L_{V^*}e^\alpha_{V^*}$ and $\sigma_K e^\alpha_K =
e^\alpha\sigma^\alpha_K$ for any open convex cone $V$ and
every closed cone
$K$. Theorem~11.5 of~\cite{Komatsu1} ensures that for an open convex
cone $V$, the boundary values of functions in $\mathcal{A}^\alpha(V)$
in the sense of ultradistributions coincide with those in the sense
of hyperfunctions. This means that $i^\alpha \mathbf
b^\alpha_V=\mathbf b_V j_V^\alpha$. It follows from these relations
and formulas~(\ref{7a}) and~(\ref{9}) that
\[
i^\alpha\mathcal F^\alpha\sigma^\alpha_{V^*} =
i^\alpha \mathbf b_V^\alpha \mathcal L^\alpha_{V^*}=
\mathbf b_V j_V^\alpha \mathcal L^\alpha_{V^*}=
\mathbf b_V\mathcal L_{V^*}e^\alpha_{V^*}=
\mathcal F \sigma_{V^*}e^\alpha_{V^*}=
\mathcal F e^\alpha\sigma^\alpha_{V^*}
\]
for any open convex cone $V$. Relation~(\ref{8}) now follows from the
inductive limit universality property.

\begin{lem}
\label{lxx8}
The canonical mapping $e^\alpha_K \colon S^{\prime
0}_\alpha(K)\to{\mathcal U}(K)$ is injective for any closed cone
$K\subset \Rset^k$.
\end{lem}

\begin{proof}
For $K=\Rset^k$, the statement follows from~(\ref{8}) because the
Fourier operator $\mathcal F^\alpha$ is an isomorphism and the
canonical mapping $i^\alpha\colon S^{\prime \alpha}_0(\Rset^k)\to
\mathcal B(\Rset^k)$ is injective by Theorem 7.5 of~\cite{Komatsu1}.
The injectivity of $e^\alpha_K$ for an arbitrary $K$ now follows from
the relation $\sigma_K e^\alpha_K = e^\alpha\sigma^\alpha_K$ and the
injectivity of the natural mapping $\sigma^\alpha_K \colon S^{\prime
0}_{\alpha}(K)\to S^{\prime 0}_{\alpha}(\Rset^k)$ ensured by
Theorem~\ref{txx0}. The lemma is proved.
\end{proof}

We end this section by establishing the connection between the
analytic wave front set (singular spectrum) of hyperfunctions and  of
their Fourier transforms.

\begin{lem}
\label{lxx10}
Let an ultrafunctional $u$ be carried by a closed convex proper cone
$K\subset \Rset^k$ and let $f=\mathcal F u$. Then the analytic wave
front set $WF_A(f)$ of the hyperfunction $f$ satisfies the relation
\[
WF_A(f)\subset \Rset^k\times (K\setminus \{0\}).
\]
\end{lem}

\begin{proof}
Theorem~9.3.3 of~\cite{Hoermander} implies that $WF_A(\mathbf b_V
\mathbf v)\subset \Rset^k\times (V^*\setminus \{0\})$ for any
connected open cone $V$ and every $\mathbf v\in\mathcal A(V)$. Hence
the assertion of the lemma follows because by definition of the
Fourier operator $\mathcal F$, we have $f=\mathbf
b_{\mathrm{int}\,K^*}\,\mathcal L_K \sigma_K\, u$.
\end{proof}

Lemma~\ref{lxx10} strengthens analogous results for tempered
distributions and ultradistributions given by Lemma~8.4.17
of~\cite{Hoermander} and Lemma~2 of~\cite{Sol3} respectively.

\section{Spaces $S^0_1(K)$ over proper cones}
\label{s4}

In this section, we show that the properties a), b), and c) listed
in Section~\ref{s2} hold also for the spaces $S^{\prime 0}_1(K)$ and
the mappings $\rho_{K',\,K}$ provided that all involved cones are
proper.
The verification of these properties
constitutes the ``functional analytic'' part of the proof of
Theorems~\ref{txx3},~\ref{txx4}, and~\ref{txx5}. In the end of
this section, we prove Theorem~\ref{txx6} describing
the Laplace transformation of ultrafunctionals carried by
proper convex closed cones.

As above, let $\dv{\cdot}{\cdot}$ be a symmetric nondegenerate
bilinear form on $\Rset^k$. For any $x,y\in\Rset^k$, we have
$|\dv{x}{y}|\leq a|x||y|$, where
\begin{equation}\label{aa}
a=\sup_{|x|,\,|y|\leq 1}|\dv{x}{y}|.
\end{equation}

\begin{lem}
\label{NewLemma}
Let $A,B>0$, $U$ be a cone in $\Rset^k$, and $W$ be a conic
neighborhood of $U$. Suppose $\eta\in \Rset^k$ is such that
$|\eta|<1/Aa$, where $a$ is given by $(\ref{aa})$. Then
the function $e^{\dv{\cdot}{\eta}}f$ belongs to
$S^0_1(U)$ for any $f\in S^{0,B}_{1,A}(W)$ and the
mapping $f\to e^{\dv{\cdot}{\eta}}f$
from $S^{0,B}_{1,A}(W)$ to $S^0_1(U)$ is continuous.
\end{lem}

\begin{proof}
Let $f\in S^{0,B}_{1,A}(W)$ and $\eta\in \Rset^k$ be such that
$|\eta|<1/Aa$. Then
\[
|e^{\dv{z}{\eta}}f(z)|\leq \|f\|_{U,A,B}\exp[-(1/A-a|\eta|)|x|+
B\delta_U(x)+B|y|],\quad z=x+iy.
\]
Therefore, $e^{\dv{\cdot}{\eta}}f\in S^{0,B}_{1,A'}(W)$,
where $A'=A/(1-Aa|\eta|)$, and the mapping
$f\to e^{\dv{\cdot}{\eta}}f$
from $S^{0,B}_{1,A}(W)$ to $S^{0,B}_{1,A'}(W)$ is continuous.
It remains to note that
the space $S^{0,B}_{1,A'}(W)$ is continuously embedded into
$S^0_1(U)$.
\end{proof}

\begin{cor}
\label{NewCor}
Let $A,B>0$ and $U$ be a cone in $\Rset^k$. Then for every
$f\in S^0_1(U)$, there is an $\varepsilon>0$ such that
$fe^{\dv{\cdot}{\eta}}\in S^{0}_{1}(U)$ for any
$\eta\in\Rset^k$ with $|\eta|<\varepsilon$.
\end{cor}

Let $U$ be a cone in $\Rset^k$, $\alpha\geq 1$, and $\eta\in
\mathrm{int}\, U^*$. We denote by $M^\eta_{\alpha,\,U}$ the
continuous mapping $f\to fe^{-\dv{\cdot}{\eta}}$ from
$S^0_\alpha(U)$ to $S^0_1(U)$.

\begin{lem}
\label{lxx4}
Let $U,U'$ be nonempty proper cones in $\Rset^k$ such that $U'\subset
U$. Then the space $S^0_1(U)$ is dense in the space $S^0_1(U')$.
\end{lem}

\begin{proof}
Fix $\alpha>1$ and
let $f\in S^0_1(U')$. By Corollary~\ref{NewCor}, there is an
$\eta\in \mathrm{int}\, U^*$ such that
$fe^{\dv{\cdot}{\eta}}\in S^{0}_{1}(U')$. This means that
$f$ belongs to the
image $\mathrm{Im}\, M^\eta_{\alpha,\,U'}$ of the mapping
$M^\eta_{\alpha,\,U'}$.
It
follows from Theorem~\ref{txx0} that $S^0_\alpha(U)$ is dense in
$S^0_\alpha(U')$. Since the image of the closure of a set under a
continuous mapping is contained in the closure of its image, we have
the inclusions $\mathrm{Im}\,
M^\eta_{\alpha,\,U'}\subset \overline{\mathrm{Im}\,
M^\eta_{\alpha,\,U}}\subset \overline{S^0_1(U)}$, where the bar
stands for closure in $S^{0}_{1}(U')$. This implies that
$f\in \overline{S^0_1(U)}$. Since $f\in S^0_1(U')$ is arbitrary,
the lemma is proved.
\end{proof}

\begin{cor}
\label{c1}
Let $K,K'$ be closed proper cones in $\Rset^k$ such that $K'\subset
K$. Then the natural mapping $\rho_{K',\,K}\colon S^{\prime
0}_1(K')\to S^{\prime 0}_1(K)$ is injective.
\end{cor}

\begin{cor}
\label{c1a}
Let $U$ be a cone in $\Rset^k$ and $U'$ be a
proper cone containing a conic neighborhood of
$U$. A functional $u\in S^{\prime 0}_1(U')$ has a continuous
extension to $S^{0}_1(U)$ if and only if $u$ can be
extended to every space $S^{0}_1(W)$, where $W$ is
a conic neighborhood of $U$ contained in $U'$.
\end{cor}

\begin{proof}
Only the sufficiency part of the statement needs
proving.
If $W\subset U'$ is a conic neighborhood of $U$, we denote
by $u_W$ the extension of $u$ to $S^0_1(W)$. By
Lemma~\ref{lxx4}, the functionals $u_W$ are uniquely
defined and are compatible with the inclusion mappings
(i.e., if $U\subset W\subset W'\subset U'$, then $u_{W'}$
is the restriction of $u_W$ to $S^0_1(W')$). As mentioned
in Section~\ref{s3}, Lemma~\ref{l0000} remains valid for
$\alpha=0$. Moreover, the union in~(\ref{OpPr1})
obviously can be taken only over conic neighborhoods of
$U$ contained in $U'$. The functionals $u_W$ therefore
determine a functional
$\tilde u\in S^{\prime 0}_1(U)$ such that
$u_W$ are restrictions of $\tilde u$ to $S^0_1(W)$.
Since $u_W$ are extensions of $u$, we conclude that
$\tilde u$ is an extension of $u$ and the corollary is
proved.
\end{proof}

\begin{lem}
\label{lxx5}
Let $K_1$ and $K_2$ be nonempty proper closed cones in $\Rset^k$ such
that $K_1\cup K_2$ is a proper cone. Then for every $f\in
S^0_1(K_1\cap K_2)$, there are $f_{1,2}\in S^0_1(K_{1,2})$ such that
$f=f_1+f_2$.
\end{lem}

\begin{proof}
Let $A,B>0$ be such that $f\in S^{0,B}_{1,A}(K_1\cap K_2)$. Fix
$\alpha>1$ and choose $\eta\in \Rset^k$ such that $|\eta|\leq 1/Aa$
and $\eta\in \mathrm{int}\, (K_1\cup K_2)^*$. Then the function
$g=fe^{\dv{\cdot}{\eta}}$ belongs to $S^0_1(K_1\cap K_2)$ and,
consequently, to $S^0_\alpha(K_1\cap K_2)$. As shown in~\cite{Sol1}
(see also Lemma~1 of~\cite{Smirnov}), there are $g_{1,2}\in
S^0_\alpha(K_{1,2})$ such that $g=g_1+g_2$. Set
$f_{1,2}=g_{1,2}e^{-\dv{\cdot}{\eta}}$. Then $f_{1,2}\in
S^0_1(K_{1,2})$ and $f=f_1+f_2$. The lemma is proved.
\end{proof}

\begin{lem}
\label{lxx6}
Let $K_1$ and $K_2$ be closed cones in $\Rset^k$. Then for every
$u\in S^{\prime 0}_1(K_1\cup K_2)$, one can find $u_{1,2}\in
S^{\prime 0}_1(K_{1,2})$ such that $u=\rho_{K_1,\,K_1\cup K_2} u_1+
\rho_{K_2,\,K_1\cup K_2}u_2$.
\end{lem}

\begin{proof}
Let $l\colon S^{0}_1(K_1\cup K_2)\to S^{0}_1(K_1)\oplus S^{0}_1(K_2)$
and $m\colon S^{0}_1(K_1)\oplus S^{0}_1(K_2) \to S^{0}_1(K_1\cap
K_2)$ be the continuous linear mappings taking $f$ to $(f,f)$ and
$(f_1,f_2)$ to $f_1-f_2$ respectively. The mapping $l$ has a closed
image because by Definition~\ref{dxx2}, we have $S^{0}_1(K_1)\cap
S^{0}_1(K_2) = S^{0}_1(K_1\cup K_2)$ and, therefore, $\mathrm{Im}\,l
=\mathrm{Ker}\, m$. In view of Lemma~\ref{lxxx0} this implies that
the space $\mathrm{Im}\,l$ is DFS.\footnote{Recall that the direct
sum of a finite family of DFS spaces and a closed subspace of a DFS
space are again DFS spaces, see~\cite{Komatsu}.} By the open mapping
theorem, the linear functional $(f,f)\to u(f)$ is continuous on
$\mathrm{Im}\,l$ and by the Hahn--Banach theorem, there exists a
continuous extension $v$ of this functional to the whole of
$S^{0}_1(K_1)\oplus S^{0}_1(K_2)$. Let $u_1$ and $u_2$ be the
restrictions of $v$ to $S^{0}_1(K_1)$ and $S^{0}_1(K_2)$
respectively. Then for any $f\in S^{0}_1(K_1\cup K_2)$, we have
$u(f)=v(f,f)=u_1(f)+u_2(f)$. This means that $u=\rho_{K_1,\,K_1\cup
K_2} u_1+ \rho_{K_2,\,K_1\cup K_2}u_2$ and the lemma is proved.
\end{proof}

\begin{lem}
\label{lxx7}
Let $K_1$ and $K_2$ be proper closed cones in $\Rset^k$ such that
$K_1\cup K_2$ is a proper cone. Let $u_{1,2}\in S^{\prime
0}_1(K_{1,2})$ be such that $\rho_{K_1,\,K_1\cup K_2} u_1=
\rho_{K_2,\,K_1\cup K_2}u_2$. Then there is a $u\in S^{\prime
0}_1(K_1\cap K_2)$ such that $u_1=\rho_{K_1\cap K_2,\,K_1}\,u$ and
$u_2=\rho_{K_1\cap K_2,\,K_2}\,u$.
\end{lem}

\begin{proof}
Let the mappings $l$ and $m$ be as in the proof of Lemma~\ref{lxx6}.
By Lemma~\ref{lxx5}, the mapping $m$ is surjective and by the open
mapping theorem, $S^{0}_1(K_1\cap K_2)$ is topologically isomorphic
to the quotient space $(S^{0}_1(K_1)\oplus
S^{0}_1(K_2))/\mathrm{Ker}\,m$. Let $v$ be the continuous linear
functional on $S^{0}_1(K_1)\oplus S^{0}_1(K_2)$ defined by the
relation $v(f_1,f_2)=u_1(f_1)-u_2(f_2)$. The condition
$\rho_{K_1,\,K_1\cup K_2} u_1= \rho_{K_2,\,K_1\cup K_2}u_2$ means
that $u_1(f)=u_2(f)$ for every $f\in S^0_1(K_1\cup K_2)$. We
therefore have $\mathrm{Ker}\,v\supset \mathrm{Im}\,l$. Since
$\mathrm{Ker}\,m = \mathrm{Im}\,l$ (see the proof of
Lemma~\ref{lxx6}), this inclusion implies the existence of a
functional $u\in S^{\prime 0}_1(K_1\cap K_2)$ such that $v=u m$. If
$f_{1,2}\in S^0_1(K_{1,2})$, then we have $u(f_1)=v(f_1,0)=u_1(f_1)$
and $u(f_2)=v(0,-f_2)=u_2(f_2)$. The lemma is proved.
\end{proof}

\begin{cor}
\label{c2}
Let $\{K_\omega\}_{\omega\in \Omega}$ be a family of closed cones in
$\Rset^k$, $K\subset \Rset^k$ be a proper closed cone such that
$K_\omega\subset K$ for every $\omega\in\Omega$, and $\tilde
K=\bigcap_{\omega\in\Omega}K_\omega$. Let $\{u_\omega\}_{\omega\in
\Omega}$ be a family of functionals such that $u_\omega\in S^{\prime
0}_1(K_\omega)$ and $\rho_{K_\omega,\,K}\,u_\omega=
\rho_{K_{\omega'},\,K}\,u_{\omega'}$ for every $\omega, \omega'\in
\Omega$. Then there is a $u\in S^{\prime 0}_1(\tilde K)$ such that
$u_\omega= \rho_{\tilde K,\,K_\omega}u$ for every $\omega\in \Omega$.
\end{cor}

\begin{proof}
If $\Omega$ is finite, then the statement follows
by induction from Lemmas~\ref{lxx7} and~\ref{lxx4}.
Now let $\Omega$ be arbitrary and
$K'\supset K$ be a proper closed cone
containing a conic neighborhood of $\tilde K$. Clearly,
the functional $v=\rho_{K_\omega,\,K'}(u_\omega)$ does not depend
on the choice of $\omega\in \Omega$.
Let $W\subset K'$ be a conic neighborhood
of $\tilde K$. By standard compactness arguments,
there is a finite family $\omega_1,\ldots,\omega_n\in\Omega$ such
that $M=\bigcap_{j=1}^n K_{\omega_j}\subset W$.
Since this corollary holds for
finite $\Omega$, we conclude that
$v$ has a continuous extension to $S^0_1(M)$ and, therefore, to
$S^0_1(W)$. Corollary~\ref{c1a} now ensures that $v$ has a
continuous extension $u$ to $S^0_1(\tilde K)$.
By Lemma~\ref{lxx4},
$\rho_{\tilde K,\,K_\omega}u$ coincides with
$u_\omega$ for any $\omega\in \Omega$
because both functionals have the same restriction to
$S^0_1(K')$. The corollary is proved.
\end{proof}

To prove Theorem~\ref{txx6}, we need the following lemma.

\begin{lem}
\label{lxx9}
Let $V$ be a convex open cone in $\Rset^k$ and $K=V^*$. Suppose a
mapping $V\ni y\to u^y\in S^{\prime 0}_1(K)$ is such that
$e^{-\dv{\cdot}{\eta'}}u^\eta=u^{\eta+\eta'}$ for any $\eta,\eta'\in
V$. Then there is a unique $u\in S^{\prime 0}_1(K)$ such that
$u^\eta=e^{-\dv{\cdot}{\eta}}u$ for any $\eta\in V$.
\end{lem}

\begin{proof}
Let $A,B>0$,
$W$ be a conic neighborhood of $K$,
and $\eta\in\Rset^k$
be such that $|\eta|\leq 1/Aa$, where $a$ is defined by
(\ref{aa}).
We denote by $L_{W,A,B}^{\eta}$ the mapping
$f\to e^{\dv{\cdot}{\eta}}f$ from $S^{0,B}_{1,A}(W)$ to $S^0_1(K)$.
Lemma~\ref{NewLemma} shows that this mapping
is well defined and continuous.
Let $\eta\in V$ be such that $|\eta|\leq 1/Aa$. We
define the continuous functional $u_{W,A,B}$ on $S^{0,B}_{1,A}(W)$ by
the relation
\[
u_{W,A,B}(f)= u^\eta(L_{W,A,B}^{\eta}\,f),\quad f\in
S^{0,B}_{1,A}(W).
\]
Although $\eta$ enters in the expression in the right-hand side,
$u_{W,A,B}$ actually does not depend on the choice of $\eta$. Indeed,
let $\eta'\in V$ be such that $|\eta'|\leq 1/Aa$. Set
$\eta''=t\eta'$, where $0<t<1$. Since $V$ is open, $\eta-\eta''\in V$
for $t$ sufficiently small, and we have
\begin{multline}
u^\eta(L_{W,A,B}^{\eta}\,f)=u^{\eta''+
(\eta-\eta'')}(e^{\dv{\cdot}{\eta}}f)=
u^{\eta''}(e^{\dv{\cdot}{\eta''}}f)=\\
=u^{\eta''+(\eta'-\eta'')}(e^{\dv{\cdot}{\eta'}}f)=
u^{\eta'}(L_{W,A,B}^{\eta'}\,f)
\nonumber
\end{multline}
for any $f\in S^{0,B}_{1,A}(W)$. Let $A'>A$, $B'>B$, and
$W'\subset W$. If $\eta\in V$ satisfies the bound $|\eta|\leq
1/A'a$, then $L_{W,A,B}^{\eta}$ is the restriction of
$L_{W',A',B'}^{\eta}$ to $S^{0,B}_{1,A}(W)$ and we have
\[
u_{W,A,B}(f) =u^\eta(L_{W,A,B}^{\eta}\,f)=
u^\eta(L_{W',A',B'}^{\eta}\,f)=
u_{W',A',B'}(f),\quad f\in S^{0,B}_{1,A}(W).
\]
Thus, the functionals $u_{W,A,B}$ are compatible with the embeddings
$S^{0,B}_{1,A}(W)\to S^{0,B'}_{1,A'}(W')$ and, therefore, determine a
functional $u\in S^{\prime 0}_1(K)$. Let $\eta,\eta'\in V$ be such
that $|\eta'|\leq 1/Aa$ and $\eta-\eta'\in V$. Fix $f\in
S^{0}_{1}(K)$ and choose $A,B>0$ and a conic neighborhood
$W$ of $K$ such that the function
$e^{-\dv{\cdot}{\eta}}f$ belongs to $S^{0,B}_{1,A}(W)$. We then
obtain
\[
(e^{-\dv{\cdot}{\eta}}u)(f)=u_{W,A,B}(e^{-\dv{\cdot}{\eta}}f)=
u^{\eta'}(e^{-\dv{\cdot}{\eta-\eta'}}f)=u^\eta(f).
\]
This means that $e^{-\dv{\cdot}{\eta}}u=u^\eta$. It remains to prove
the uniqueness of $u$. Suppose $u'\in S^{\prime 0}_1(K)$ is such that
$e^{-\dv{\cdot}{\eta}}u'=u^\eta$ for any $\eta\in V$. Then $v=u'-u$
satisfies the relation $e^{-\dv{\cdot}{\eta}}v=0$ for any $\eta\in
V$. Let $f\in S^0_1(K)$. By Corollary~\ref{NewCor}, there is an
$\eta\in V$ such that $e^{\dv{\cdot}{\eta}}f\in S^0_1(K)$. We
therefore have
$v(f)=(e^{-\dv{\cdot}{\eta}}v)(e^{\dv{\cdot}{\eta}}f)=0$. Thus, $v=0$
and the lemma is proved.
\end{proof}

\begin{proof}[Proof of Theorem~$\ref{txx6}$.]
As in the preceding section, we identify ${\mathcal U}(K)$ with
$S^{\prime 0}_1(K)$. Fix $\alpha>1$. For $u\in S^{\prime 0}_1(K)$ and
$\eta\in V$, we define the functional $v^\eta\in S^{\prime
0}_\alpha(K)$ by the relation $v^\eta(f)=u(e^{-\dv{\cdot}{\eta}}f)$,
$f\in S^0_\alpha(K)$. We then have
\begin{equation}\label{10}
(\mathcal L_K u)(\zeta+i\eta)=v^\eta(e^{i\dv{\cdot}{\zeta}})=
(\mathcal L^\alpha_K v^\eta)(\zeta), \quad \zeta\in T^V,
\end{equation}
and in view of Theorem~\ref{txx8} the function $(\mathcal L_K
u)(\cdot+i\eta)$ is analytic in $T^V$. Since $\eta\in V$ is
arbitrary, this implies that $\mathcal L_K u$ is analytic in $T^V$.
If $\mathcal L_K u=0$, then by~(\ref{10}) we have $\mathcal
L^\alpha_K v^\eta=0$ for any $\eta\in V$. This implies that
$v^\eta=0$ for any $\eta \in V$ because the Laplace transformation
$\mathcal L^\alpha_K$ is injective by Theorem~\ref{txx8}. Denoting by
$u^\eta$ the restriction of $v^\eta$ to $S^0_1(K)$ and applying the
uniqueness part of Lemma~\ref{lxx9}, we conclude that $u=0$. Thus,
the operator $\mathcal L_K$ is injective. The mapping $u\to v^\eta$
from $S^{\prime 0}_1(K)$ to $S^{\prime 0}_\alpha(K)$ is continuous
for any $\eta\in V$ being the dual mapping of the continuous mapping
$f\to e^{-\dv{\cdot}{\eta}}f$.
It therefore follows from~(\ref{10}) and Theorem~\ref{txx8} that the
mapping $u\to (\mathcal L_K u)(\cdot+i\eta)$ is continuous as a
mapping from $S^{\prime 0}_1(K)$ to $\mathcal A^\alpha(V)$ and,
consequently, as a mapping from $S^{\prime 0}_1(K)$ to $\mathcal
A(V)$. This implies the continuity of $\mathcal L_K$ because for
every compact set $K\subset T^V$, one can find an $\eta\in V$ such
that $K-\eta\subset T^V$. We now prove the surjectivity of $\mathcal
L_K$. Let $\mathbf v\in\mathcal A(V)$. Clearly, $\mathbf
v(\cdot+i\eta)\in \mathcal A^\alpha(V)$ for any $\eta\in V$. We set
$v^\eta=(\mathcal L^\alpha_K)^{-1}\mathbf v(\cdot+i\eta)$ and denote
by $u^\eta$ the restrictions of $v^\eta$ to $S^0_1(K)$. For every
$\eta\in V$ and $\zeta\in T^V$, we have
\[
(\mathcal L_K u^\eta)(\zeta)= u^{\eta}(e^{i\dv{\cdot}{\zeta}})=
v^\eta(e^{i\dv{\cdot}{\zeta}})=
(\mathcal L^\alpha_K v^\eta)(\zeta)=\mathbf v(\zeta+i\eta).
\]
Hence it follows that
\[
\mathcal L_K (e^{-\dv{\cdot}{\eta'}} u^\eta)=
(\mathcal L_K u^\eta)(\cdot +i\eta')=
\mathbf v(\cdot +i(\eta+\eta'))=
\mathcal L_K u^{\eta+\eta'}, \quad \eta,\eta'\in V,
\]
and in view of the injectivity of $\mathcal L_K$ we have
$e^{-\dv{\cdot}{\eta'}} u^\eta=u^{\eta+\eta'}$. By Lemma~\ref{lxx9},
there is a $u\in S^{\prime 0}_1(K)$ such that $u^\eta=
e^{-\dv{\cdot}{\eta}}u$ for any $\eta\in V$. Fix $\zeta=\xi+i\eta\in
T^V$ and choose $\eta'\in V$ such that $\eta-\eta'\in V$. Then we
have
\[
(\mathcal L_K u)(\zeta)=
(e^{-\dv{\cdot}{\eta'}}u)(e^{i\dv{\cdot}{\xi+i(\eta-\eta')}})=
(\mathcal L_K u^{\eta'})(\xi+i(\eta-\eta'))=\mathbf v(\zeta).
\]
Thus, $\mathcal L_K u=\mathbf v$ and, consequently, $\mathcal L_K$ is
a continuous one-to-one mapping. Since both $S^{\prime 0}_1(K)$ and
$\mathcal A(V)$ are Fr\'echet spaces, the continuity of the inverse
operator $\mathcal L_K^{-1}$ is ensured by the open mapping theorem.
Theorem~$\ref{txx6}$ is proved.
\end{proof}

\section{Localizable inductive systems}
\label{s5}

The results of the preceding section show that the localization
properties described by Theorems~\ref{txx3},~\ref{txx4},
and~\ref{txx5} hold for ultrafunctionals carried by proper closed
cones. To prove these theorems in their full volume, we have to show
that the properties of the inductive system $S$ formed by the spaces
$S^{\prime 0}_1(K)$ over proper closed cones are inherited by the
inductive system $\mathcal U$ formed by the spaces $\mathcal U(K)$
over arbitrary closed cones. We shall obtain the desired localization
properties of $\mathcal U$ as a consequence of a more general
algebraic construction formulated in terms of \emph{(pre)localizable
inductive systems} introduced by Definition~\ref{dxx6} below. In
contrast to Section~\ref{s4}, all considerations in this section are
purely algebraic.

Recall that a partially ordered set $A$ is called a lattice if every
two-element subset $\{\alpha_1,\alpha_2\}$ of the set $A$ has a
supremum $\alpha_1\vee \alpha_2$ and an infimum $\alpha_1\wedge
\alpha_2$. A lattice $A$ is called distributive if
$\alpha_1\wedge(\alpha_2\vee\alpha_3)=
(\alpha_1\wedge\alpha_2)\vee(\alpha_1\wedge\alpha_3)$
for any $\alpha_1,\alpha_2,\alpha_3\in A$.

\begin{defn}
\label{dxx4}
Let $A$ be a partially ordered set. We say that $A$ is a
quasi-lattice if every two-element subset of $A$ has an infimum and
every bounded above two-element subset of $A$ has a supremum. We say
that a quasi-lattice $A$ is distributive if
$\alpha_1\wedge(\alpha_2\vee\alpha_3)=
(\alpha_1\wedge\alpha_2)\vee(\alpha_1\wedge\alpha_3)$
for every bounded above pair $\alpha_2,\alpha_3\in A$ and every
$\alpha_1\in A$.
\end{defn}

Clearly, every (distributive) lattice is a (distributive)
quasi-lattice. If $A$ is a distributive lattice, then we have
\begin{multline}
(\alpha_1\vee\alpha_2)\wedge(\alpha_1\vee\alpha_3)=
((\alpha_1\vee\alpha_2)\wedge\alpha_1)\vee
((\alpha_1\vee\alpha_2)\wedge\alpha_3)=\\
=\alpha_1\vee((\alpha_1\wedge\alpha_3)\vee
(\alpha_2\wedge\alpha_3))=\alpha_1\vee (\alpha_2\wedge\alpha_3)
\nonumber
\end{multline}
for any $\alpha_1,\alpha_2,\alpha_3\in A$. By induction, it follows
that
\begin{equation}\label{0}
\inf_{\omega\in \Omega}(\alpha\vee\alpha_\omega)=
\alpha\vee \inf_{\omega\in \Omega}\alpha_\omega
\end{equation}
for any $\alpha\in A$ and every nonempty finite family
$\{\alpha_\omega\}_{\omega\in \Omega}$ of elements of $A$.

\begin{defn}
\label{dxx5}
We call a lattice $A$ infinitely distributive if every nonempty
subset of $A$ has an infimum and the condition~(\ref{0}) is satisfied
for an arbitrary (not necessarily finite) family
$\{\alpha_\omega\}_{\omega\in \Omega}$ of elements of $A$.
\end{defn}

Note that a distributive lattice may be not infinitely distributive
even if all its subsets have an infimum (see, e.g.,~\cite{Graetzer},
Section~II.4, exercises~17 and~18).

We call a nondecreasing mapping $\lambda$ from a quasi-lattice $A$
into a quasi-lattice $B$ a morphism of quasi-lattices if
$\lambda(\alpha_1\wedge\alpha_2)=
\lambda(\alpha_1)\wedge\lambda(\alpha_2)$ for any
$\alpha_1,\alpha_2\in A$ and $\lambda(\alpha_1\vee\alpha_2)=
\lambda(\alpha_1)\vee\lambda(\alpha_2)$ for every bounded above pair
$\alpha_1,\alpha_2\in A$.

In the rest of this section, we study abstract inductive systems of
vector spaces indexed by (quasi-)lattices and systematically use the
corresponding notation introduced in the end of Section~\ref{s2}.

\begin{defn}
\label{dxx6}
An inductive system $\mathcal X$ of vector spaces over a
quasi-lattice $A$ is called to be prelocalizable if the following
conditions are satisfied:
\begin{enumerate}
\item[I] The mappings $\rho^{\mathcal X}_{\alpha\alpha'}$ are
injective for any $\alpha,\alpha'\in A$, $\alpha\leq\alpha'$.
\item[II] If a pair $\alpha_1,\alpha_2\in A$ is bounded above and
$x\in \mathcal X(\alpha_1\vee\alpha_2)$, then there are $x_{1,\,2}\in
\mathcal X(\alpha_{1,\,2})$ such that $x=\rho^{\mathcal
X}_{\alpha_1,\,\alpha_1\vee\alpha_2}(x_1)+ \rho^{\mathcal
X}_{\alpha_2,\,\alpha_1\vee\alpha_2}(x_2)$.
\item[III] If a pair $\alpha_1,\alpha_2\in A$ is bounded above by an
element $\alpha\in A$, $x_{1,\,2}\in\mathcal X(\alpha_{1,\,2})$, and
$\rho^{\mathcal X}_{\alpha_1,\,\alpha}(x_1)= \rho^{\mathcal
X}_{\alpha_2,\,\alpha}(x_2)$, then there is an $x\in \mathcal
X(\alpha_1\wedge\alpha_2)$ such that $x_1=\rho^{\mathcal
X}_{\alpha_1\wedge\alpha_2,\,\alpha_1}(x)$ and $x_2=\rho^{\mathcal
X}_{\alpha_1\wedge\alpha_2,\,\alpha_2}(x)$.
\end{enumerate}
We say that the inductive system $\mathcal X$ is localizable if every
nonempty subset of $A$ has an infimum and instead of III the
following stronger condition is satisfied:
\begin{enumerate}
\item[III$'$] Let $\{\alpha_\omega\}_{\omega\in \Omega}$ be a
nonempty family of elements of $A$ bounded above by an $\alpha\in A$,
and let a family $\{x_\omega\}_{\omega\in \Omega}$ be such that
$x_\omega\in\mathcal X(\alpha_\omega)$ and $\rho^{\mathcal
X}_{\alpha_\omega,\,\alpha}(x_\omega)= \rho^{\mathcal
X}_{\alpha_{\omega'},\,\alpha}(x_{\omega'})$ for any $\omega,
\omega'\in \Omega$. Then there is an $x\in \mathcal X(\tilde \alpha)$
($\tilde \alpha=\inf_{\omega\in \Omega} \alpha_\omega$) such that
$x_\omega= \rho^{\mathcal X}_{\tilde\alpha,\,\alpha_\omega}(x)$ for
any $\omega\in \Omega$.
\end{enumerate}
\end{defn}

Let $M$ be a closed cone in $\Rset^k$. The (ordered by inclusion) set
$\mathcal P(M)$ of all proper closed cones contained in $M$ is a
distributive quasi-lattice, while the set $\mathcal K(M)$ of all
closed cones contained in $M$ is an infinitely distributive lattice.
As shown by the properties a)--c) listed in Section~\ref{s2}, the
inductive system over $\mathcal K(\Rset^k)$ formed by the spaces
$S^{\prime 0}_\alpha(K)$ ($\alpha>1$) is prelocalizable (in fact, it
is even localizable, see the paragraph following the formulation of
Theorem~\ref{txx2}).

A subset $I$ of a quasi-lattice $A$ will be called $\wedge$-closed if
$\alpha_1\wedge\alpha_2\in I$ for any $\alpha_1,\alpha_2\in I$. If
$I$ is a finite subset of a quasi-lattice $A$, then one can find a
finite $\wedge$-closed set $I'\subset A$ containing $I$ (for
instance, the set consisting of infima of all subsets of $I$ can be
taken as $I'$).

\begin{lem}
\label{l00}
Let $A$ be a distributive lattice, $\alpha\in A$, and $\mathcal X$ be
a prelocalizable inductive system over $\mathcal X$. Suppose $I$ is a
$\wedge$-closed subset of $A$ such that $\alpha'\leq\alpha$ for any
$\alpha'\in I$ and $\alpha=\alpha_1\vee\ldots\vee\alpha_n$ for some
$\alpha_1,\ldots,\alpha_n\in I$. Then the space $\mathcal X(\alpha)$
is canonically isomorphic to $\varinjlim_{\alpha'\in I}\mathcal
X(\alpha')$.
\end{lem}

The proof of this lemma is completely analogous to the algebraic part
of the proof of Lemma~\ref{lxx0} and is omitted. The following result
is an immediate consequence of Corollary~\ref{c1},
Lemma~\ref{lxx6}, and Corollary~\ref{c2}.

\begin{lem}
\label{l0}
The inductive system $S$ over $\mathcal P(\Rset^k)$ formed by the
spaces $S^{\prime 0}_1(K)$ is localizable.
\end{lem}

Theorems~\ref{txx3},~\ref{txx4}, and~\ref{txx5} can be reformulated
in terms of localizable inductive systems as follows.

\begin{thm}
\label{t0}
The inductive system ${\mathcal U}$ over $\mathcal K(\Rset^k)$ formed
by the spaces ${\mathcal U}(K)$ is localizable.
\end{thm}

Let $\mathcal X$ be an inductive system over a partially ordered set
$A$. For every $I\subset A$, we define the inductive system $\mathcal
X^I$ over $I$ setting $\mathcal X^I(\alpha)=\mathcal X(\alpha)$ and
$\rho^{\mathcal X^I}_{\alpha\alpha'}=\rho^{\mathcal
X}_{\alpha\alpha'}$ for $\alpha,\alpha'\in I$, $\alpha\leq \alpha'$
(i.e., $\mathcal X^I$ is the ``restriction'' of $\mathcal X$ to $I$).
If $I\subset J\subset A$, then there are canonical mappings
$\tau^{\mathcal X}_{I,\,J}\colon \varinjlim \mathcal X^I\to
\varinjlim \mathcal X^J$ satisfying the relation $\tau^{\mathcal
X}_{I,\,J}\rho_\alpha^{\mathcal X^I}= \rho_\alpha^{\mathcal X^J}$ for
any $\alpha\in I$. Let $\lambda$ be a nondecreasing mapping from $A$
to a partially ordered set $B$. With every $\beta \in B$ we associate
the set $A_\beta=\{\alpha\in A\,|\,\lambda(\alpha)\leq \beta\}$ and
define the inductive system $\lambda(\mathcal X)$ over $B$ setting
$\lambda(\mathcal X)(\beta)=\varinjlim \mathcal X^{A_\beta}$ and
$\rho^{\lambda(\mathcal X)}_{\beta\beta'}= \tau^{\mathcal
X}_{A_\beta,\,A_{\beta'}}$ for $\beta,\beta'\in B$,
$\beta\leq\beta'$.

The inclusion mapping $\theta\colon\mathcal P(\Rset^k)\to \mathcal
K(\Rset^k)$ is clearly a morphism of quasi-lattices. By definition of
the inductive system ${\mathcal U}$, we have ${\mathcal U}=
\theta(S)$ and, therefore, Theorem~\ref{t0} follows from the
following more general statement.

\begin{thm}
\label{t1}
Let $A$ be a distributive quasi-lattice, $B$ be a distributive
lattice and $\lambda\colon A\to B$ be an injective quasi-lattice
morphism such that every element $\beta\in B$ is representable in the
form $\beta=\lambda(\alpha_1)\vee\ldots\vee\lambda(\alpha_n)$, where
$\alpha_1,\ldots,\alpha_n\in A$. If $\mathcal X$ is a prelocalizable
inductive system of vector spaces over $A$, then $\lambda(\mathcal
X)$ is a prelocalizable inductive system of vector spaces over $B$.
If $\mathcal X$ is a localizable inductive system over $A$ and the
lattice $B$ is infinitely distributive, then the inductive system
$\lambda(\mathcal X)$ is localizable.
\end{thm}

The proof of Theorem~\ref{t1} is given in Appendix~\ref{app1}.

\begin{rem}
Under the conditions of Theorem~\ref{t1}, for every quasi-lattice
morphism $\varphi$ from $A$ to a distributive lattice $L$, there is a
unique lattice morphism $\psi\colon B\to L$ such that
$\varphi=\psi\lambda$. This means that $B$ is the \emph{free
distributive lattice} over the partially ordered set $A$
(see~\cite{Graetzer}, Sec.~I.5, Definition~2).
\end{rem}

Let $K_1,\ldots,K_n$ be convex closed proper cones in $\Rset^k$ such
that $\bigcup_{i=1}^n K_i=\Rset^k$ and let $I$ be the set consisting
of all intersections of the cones $K_1,\ldots,K_n$. It follows from
Theorems~\ref{txx6} and~\ref{t0} and Lemma~\ref{l00} that $\mathcal
U(\Rset^k)$ is canonically isomorphic to the space $\varinjlim_{K\in
I} \mathcal A(\mathrm{int}\,K^*)$. In the next section, we shall
establish the bijectivity of the Fourier transformation by proving
that for some choice of the cones $K_i$, the latter space is
isomorphic to $\mathcal B(\Rset^k)$. To this end, we shall need to
pass from the above inductive limit representation of $\mathcal
U(\Rset^k)$ to another representation similar to that given by
Martineau's edge of the wedge theorem for hyperfunctions. We conclude
this section by describing the corresponding procedure in terms of
abstract inductive systems.

Recall~\cite{Graetzer} that a partially ordered set $A$ is called a
lower semilattice if every two-element subset $\{\alpha_1,\alpha_2\}$
of the set $A$ has an infimum $\alpha_1\wedge \alpha_2$. Recall also
that a subset $I$ of a partially ordered set $A$ is called cofinal in
$A$ if every element of $A$ is majorized by an element of $I$.

\begin{lem}
\label{l5}
Let $T$ be a set, $A$ be a lower semilattice, $\mathcal X$ be an
inductive system over $A$, and $\lambda$ be a mapping from $T$ to $A$
such that $\lambda(T)$ is cofinal in $A$. Let $\EuScript N$ be the
subspace of $\oplus_{\tau\in T} \mathcal X(\lambda(\tau))$ spanned by
all vectors of the form
\begin{equation}\label{11}
j_\tau \rho^{\mathcal X}_{\lambda(\tau)\wedge
\lambda(\tau'),\,\lambda(\tau)} x-
j_{\tau'} \rho^{\mathcal X}_{\lambda(\tau)\wedge
\lambda(\tau'),\,\lambda(\tau')} x,
\end{equation}
where
$\tau,\tau'\in T$, $x\in \mathcal
X(\lambda(\tau)\wedge\lambda(\tau'))$, and $j_\tau$ is the canonical
embedding of $\mathcal X(\lambda(\tau))$ into $\oplus_{\tau\in T}
\mathcal X(\lambda(\tau))$. Then we have a natural isomorphism
\[
\oplus_{\tau\in T} \mathcal X(\lambda(\tau))/\EuScript N\simeq
\varinjlim \mathcal X.
\]
\end{lem}

The proof of Lemma~\ref{l5} is given in Appendix~\ref{app2}.

\begin{cor}
\label{c3}
Let $A$ be a finite lower semilattice and $\mathcal X$ be an
inductive system over $A$. Let $\alpha_1,\ldots,\alpha_n$ be a family
of elements of $A$ containing all maximal elements of $A$ and
$\EuScript N$ be the subspace of $\oplus_{i=1}^n \mathcal
X(\alpha_i)$ consisting of the vectors $(x_1,\ldots,x_n)$ whose
components are representable in the form
\begin{equation}
\label{edge}
x_i=\sum_{l=1}^n j_{l}
\rho^{\mathcal X}_{\alpha_i\wedge\alpha_l,\,\alpha_l}
x_{il},\quad i=1,\ldots,n,
\end{equation}
where $x_{il}\in \mathcal X(\alpha_i\wedge \alpha_l)$,
$x_{il}=-x_{li}$, and $j_i$ is the canonical embedding of $\mathcal
X(\alpha_i)$ into $\oplus_{i=1}^n \mathcal X(\alpha_i)$. Then we have
a natural isomorphism
\[
\oplus_{i=1}^n \mathcal X(\alpha_i)/\EuScript N\simeq
\varinjlim \mathcal X.
\]
\end{cor}

\begin{proof}
Obviously, a subset $I$ of a finite partially ordered set $A$ is
cofinal in $A$ if and only if it contains all maximal elements of
$A$, and in view of Lemma~\ref{l5} it suffices to show that
$\EuScript N$ coincides with the subspace $\EuScript N'$ of
$\oplus_{i=1}^n \mathcal X(\alpha_i)$ spanned by the vectors of the
form $j_{i}\rho^{\mathcal X}_{\alpha_i\wedge\alpha_l,\,\alpha_i} y-
j_{l}\rho^{\mathcal X}_{\alpha_i\wedge\alpha_l,\,\alpha_l} y$ with
$y\in \mathcal X(\alpha_i\wedge \alpha_l)$. Let $x\in \EuScript N'$.
Then we have
\[
x=\sum_{i,l=1}^n (j_{i}\rho^{\mathcal X}_{\alpha_i\wedge
\alpha_l,\,\alpha_i} y_{il}-
j_{l}\rho^{\mathcal X}_{\alpha_i\wedge\alpha_l,\,\alpha_l} y_{il})=
\sum_{i=1}^n j_i \sum_{l=1}^n \rho^{\mathcal X}_{\alpha_i\wedge
\alpha_l,\,\alpha_i} (y_{il}-y_{li}),
\]
where $y_{il}\in \mathcal X(\alpha_i\wedge \alpha_l)$. Setting
$x_{il}=y_{il}-y_{li}$, we see that the components of $x$ have the
form~(\ref{edge}) and, therefore, $x\in \EuScript N$. Conversely, let
$x$ be the element of $\EuScript N$ whose components have the
form~(\ref{edge}). Then in view of the antisymmetry of $x_{il}$ we
have
\[
x=\sum_{i=1}^n j_i \sum_{l=1}^n \rho^{\mathcal X}_{\alpha_i\wedge
\alpha_l,\,\alpha_i} x_{il}=
\frac{1}{2}\sum_{i,l=1}^n (j_{i}\rho^{\mathcal X}_{\alpha_i\wedge
\alpha_l,\,\alpha_i} x_{il}-
j_{l}\rho^{\mathcal X}_{\alpha_i\wedge\alpha_l,\,\alpha_l} x_{il})
\]
and, therefore, $x\in \EuScript N'$. Thus, $\EuScript N=\EuScript N'$
and the corollary is proved.
\end{proof}

\section{Bijectivity of Fourier transformation}
\label{s6}

In this section, we give the proof of Theorem~\ref{txx7}.

We first consider the one-dimensional case, when the spaces of
hyperfunctions and ultrafunctionals have very simple structure and
the bijectivity of $\mathcal F$ can be derived immediately from
Theorem~\ref{txx6} without any reference to algebraic constructions
of the preceding section. Let $H(V)$ denote the space of functions
holomorphic in an open set $V\subset \Cset$. According to Sato's
definition, hyperfunctions on an open set $O\subset \Rset$ are the
elements of the quotient space $H(V\setminus O)/H(V)$, where $V$ is
an open set in $\Cset$ containing $O$ as a relatively closed subset
and $H(V)$ is assumed to be embedded in $H(V\setminus O)$ via the
restriction mapping. It is important that all such quotient spaces
are naturally isomorphic to each other and, therefore, this
definition actually does not depend on the choice of $V$ (see, e.g.,
Section~2 of~\cite{Komatsu2} or Section 3.1 of~\cite{Morimoto}). In
particular, we can set $\mathcal
B(\Rset)=H(\Cset\setminus\Rset)/H(\Cset)$. For $\mathbf v\in
H(\Cset\setminus\Rset)$, we denote by $[\mathbf v]$ the corresponding
element of $\mathcal B(\Rset)$. Let the operators $j_\pm\colon
H(\Cset_\pm)\to H(\Cset\setminus\Rset)$ be defined by the relations
\[
(j_\pm\mathbf v_\pm)(\zeta)=\left\{
\begin{matrix}
\mathbf v_\pm(\zeta) & \mbox{for }\zeta\in \Cset_\pm\\
0 & \mbox{for }\zeta\in \Cset_\mp
\end{matrix}
\right., \quad \mathbf v_\pm\in H(\Cset_\pm).
\]
The boundary value operators $\mathbf b_{\Rset_\pm}\colon \mathcal
A(\Rset_\pm)\to\mathcal B(\Rset)$ and $\mathbf b_{\Rset}\colon
\mathcal A(\Rset)\to\mathcal B(\Rset)$ are given by
\begin{equation}\label{bound}
\mathbf b_{\Rset_\pm}\mathbf v_\pm = \pm [j_\pm \mathbf v_\pm], \quad
\mathbf b_{\Rset}\mathbf v = \mathbf b_{\Rset_+}(\mathbf v|_{\Cset_+})
=\mathbf b_{\Rset_-}(\mathbf v|_{\Cset_-}),
\end{equation}
where
$\mathbf v_\pm\in \mathcal A(\Rset_\pm)$,
$\mathbf v\in\mathcal A(\Rset)$, and
$\mathbf v|_{\Cset_\pm}$ are the restrictions of $\mathbf v$ to
$\Cset_\pm$ (note that $\mathcal A(\Rset_\pm)=H(\Cset_\pm)$ and
$\mathcal A(\Rset)=H(\Cset)$). By Theorem~\ref{txx6} and the
definition of $\mathcal U(\Rset)$, the Laplace operators $\mathcal
L_{K}$ determine an isomorphic mapping $\mathcal L\colon {\mathcal
U}(\Rset)\to \varinjlim_{K\in \mathcal P(\Rset)} \mathcal
A(\mathrm{int}\,K^*)$. The set $\mathcal P(\Rset)$ contains only
three elements: $\bar\Rset_+$, $\bar\Rset_-$, and $\{0\}$. By
definition of the inductive limit, we have
\[
\varinjlim_{K\in \mathcal P(\Rset)} \mathcal A(\mathrm{int}\,K^*)=
(\mathcal A(\Rset)\oplus \mathcal A(\Rset_+)\oplus
\mathcal A(\Rset_-))/N,
\]
where $N$ is the subspace of $\mathcal A(\Rset)\oplus \mathcal
A(\Rset_+)\oplus \mathcal A(\Rset_-)$ spanned by the vectors of the
form $(\mathbf v,-\mathbf v|_{\Cset_+},0)$ and $(\mathbf v,0,-\mathbf
v|_{\Cset_-})$ with $\mathbf v\in \mathcal A(\Rset)$. Let the mapping
$\tilde s\colon \mathcal A(\Rset)\oplus \mathcal A(\Rset_+)\oplus
\mathcal A(\Rset_-)\to \mathcal B(\Rset)$ be defined by the relation
\begin{equation}\label{bij}
\tilde s(\mathbf v,\mathbf v_+,\mathbf v_-)=
\mathbf b_\Rset \mathbf v + \mathbf b_{\Rset_+} \mathbf v_+ +
\mathbf b_{\Rset_-} \mathbf v_-
\end{equation}
and let $s\colon \varinjlim_{K\in \mathcal P(\Rset)} \mathcal
A(\mathrm{int}\,K^*)\to \mathcal B(\Rset)$ be the mapping induced by
$\tilde s$. By definition of the Fourier operator $\mathcal F$, we
have $\mathcal F= s\mathcal L$. Thus, to prove the bijectivity of
$\mathcal F$, it suffices to show that $s$ is one-to-one. In other
words, we have to show that $(\mathbf v,\mathbf v_+,\mathbf v_-)\in
N$ if $\tilde s(\mathbf v,\mathbf v_+,\mathbf v_-)=0$. In view
of~(\ref{bound}) and~(\ref{bij}), the last condition means that
$[j_+(\mathbf v|_{\Cset_+})]+[j_+\mathbf v_+]-[j_-\mathbf v_-]=0$. In
other words, there is a $\mathbf u\in \mathcal A(\Rset)$ such that
$\mathbf v|_{\Cset_+}+\mathbf v_+=\mathbf u|_{\Cset_+}$ and $-\mathbf
v_-=\mathbf u|_{\Cset_-}$. This implies that
\[
(\mathbf v,\mathbf v_+,\mathbf v_-)=
(\mathbf v-\mathbf u,-(\mathbf v-\mathbf u)|_{\Cset_+},0)+
(\mathbf u, 0, -\mathbf u|_{\Cset_-}).
\]
Thus, $(\mathbf v,\mathbf v_+,\mathbf v_-)\in N$ and
Theorem~\ref{txx7} is proved for the case $k=1$.

Let us now consider the general case. With every set
$x^1,\ldots,x^{l}$ of vectors in $\Rset^k$ we associate the cone
$K(x^1,\ldots,x^l)=\{\,x\in \Rset^k\,|\,x=t_1 x^1+\ldots+t_l x^l,\,
t_i\geq 0\,\}$. Let $x^1,\ldots, x^{k+1}$ be vectors in $\Rset^k$
such that $K(x^1,\ldots, x^{k+1})=\Rset^k$. For $i,j=1,\ldots,k+1$,
$i\ne j$, we set
\begin{align}
E_i &= \{\,\xi\in \Rset^k\,|\, \dv{\xi}{x_i}\geq 0\,\}, \nonumber\\
K_i &= K(x^1,\ldots,\hat x^i,\ldots x^{k+1}),\quad K_{ij}=
K(x^1,\ldots,\hat x^i,\ldots,\hat x^j,\ldots x^{k+1}),
\nonumber\\
\Gamma_i &= E_1\cap\ldots\cap\hat E_i\cap\ldots \cap E_{k+1},\quad
V_{ij}=E_1\cap\ldots\cap\hat E_i\cap\ldots\cap\hat E_j\cap\ldots
\cap E_{k+1},\nonumber
\end{align}
where $\dv{\cdot}{\cdot}$ is the symmetric nondegenerate bilinear
form on $\Rset^k$ entering in the definitions of the Fourier and
Laplace transformations and the hat means that an element is omitted.
It is easy to see that $\bigcup_{i=1}^{k+1} K_i=\Rset^k$ and $K_i\cap
K_j= K_{ij}$. Furthermore, we have $\Gamma_i = \mathrm{int}\, K_i^*$
and $V_{ij} = \mathrm{int}\, K_{ij}^*$. Let $I$ denote the finite set
consisting of all possible intersections of the cones
$K_1,\ldots,K_{k+1}$. By Lemma~\ref{l00} and Theorem~\ref{t0}, there
is a natural isomorphism $l\colon {\mathcal U}(\Rset^k)\to
\varinjlim_{K\in I} {\mathcal U}(K)$. By Theorem~\ref{txx6}, the
Laplace operators $\mathcal L_{K}$ determine an isomorphic mapping
$\mathcal L\colon \varinjlim_{K\in I} {\mathcal U}(K)\to
\varinjlim_{K\in I} \mathcal A(\mathrm{int}\,K^*)$. Let $\EuScript N$
be the subspace of $\oplus_{i=1}^{k+1} \mathcal A(\Gamma_i)$
consisting of the elements $(\mathbf v_1,\ldots,\mathbf v_{k+1})$
such that
\[
\mathbf v_i = \sum_{j=1}^{k+1} \mathbf v_{ij},
\]
where $\mathbf v_{ij}=-\mathbf v_{ji}$ belong to $\mathcal
A(V_{ij})$. By Corollary~\ref{c3}, we have a natural
isomorphism $m \colon \varinjlim_{K\in I} \mathcal
A(\mathrm{int}\,K^*)\to \oplus_{i=1}^{k+1} \mathcal
A(\Gamma_i)/\EuScript N$. Let $\tilde b$ be the mapping from
$\oplus_{i=1}^{k+1} \mathcal A(\Gamma_i)$ to $\mathcal B(\Rset^k)$
defined by the relation
\[
\tilde b(\mathbf v_1,\ldots,\mathbf v_{k+1})=
\sum_{i=1}^{k+1} \mathbf b_{\Gamma_i} \mathbf v_i,
\]
where $\mathbf b_{\Gamma_i}$ are the boundary value operators.
Obviously, we have the inclusion $\EuScript N\subset \ker \tilde b$
and, therefore, $\tilde b$ determines a mapping $b \colon
\oplus_{i=1}^{k+1} \mathcal A(\Gamma_i)/\EuScript N\to \mathcal
B(\Rset^k)$. From the definition of the Fourier operator $\mathcal
F$, it easily follows that $\mathcal F = b m \mathcal L l$. Thus, it
suffices to establish that $b$ is a one-to-one mapping. Let $\tilde
B$ be the mapping from $\oplus_{i=1}^{k+1} \mathcal A(\Gamma_i)$ to
$\mathcal B(\Rset^k)$ defined by the relation
\[
\tilde B(\mathbf v_1,\ldots,\mathbf v_{k+1})=
\sum_{i=1}^{k+1}(-1)^i \,\mathbf b_{\Gamma_i} \mathbf v_i.
\]
Let the mapping $\delta \colon \oplus_{i<j} \mathcal A(V_{ij})\to
\oplus_i \mathcal A(\Gamma_i)$ be defined by
\[
(\delta\mathbf v)_j=\sum_{1\leq i<j}(-1)^i \mathbf v_{ij}+
\sum_{j< i\leq k+1}(-1)^{i+1} \mathbf v_{ji},\quad \mathbf v=
\{\mathbf v_{ij}\}_{i<j}.
\]
It is easy to see that
$\mathrm{Im}\,\delta\subset\mathrm{Ker}\,\tilde B$ and, consequently,
$\tilde B$ determines a mapping $B\colon \oplus_{i=1}^{k+1} \mathcal
A(\Gamma_i)/\mathrm{Im}\,\delta \to \mathcal B(\Rset^k)$. As shown
in~\cite{Morimoto} (see formula~2.5 of Chapter~7 and Corollary~7.4.6)
this mapping is one-to-one. Let $\tau$ be the isomorphic mapping from
$\oplus_{i=1}^{k+1} \mathcal A(\Gamma_i)$ onto itself defined by
\[
(\tau\mathbf v)_i = (-1)^i\mathbf v_i,\quad \mathbf v=
(\mathbf v_1,\ldots,\mathbf v_{k+1}).
\]
Then we have $\tilde b=\tau \tilde B$ and
$\mathrm{Im}\,\delta=\tau(\EuScript N)$. Therefore, the bijectivity
of $B$ implies that of $b$. Theorem~\ref{txx7} is proved.

\begin{rem}
Let
\[
\mathfrak V = (\Cset^k,T^{E_1},\ldots, T^{E_{k+1}}),\quad
\mathfrak V' = (T^{E_1},\ldots, T^{E_{k+1}}),
\]
where $T^{E_j}=\Rset^k +i E_j$. The above isomorphism $B\colon
\oplus_{i=1}^{k+1} \mathcal A(\Gamma_i)/\mathrm{Im}\,\delta \to
\mathcal B(\Rset^k)$ gives the \v{C}ech cohomology representation of
$\mathcal B(\Rset^k)$ if $(\mathfrak V,\mathfrak V')$ is used as a
relative Stein open covering of $(\Cset^k,\Cset^k\setminus \Rset^k)$
(see details in~\cite{Morimoto}, Sec.~7.2).
\end{rem}

\section{Conclusion}
\label{s7}

The obtained results suggest the way of constructing
``nontrivializations'' of some seemingly trivial generalized function
spaces. We conclude this paper by indicating some possible results of
this type in the framework of the Gurevich spaces $W^\Omega_M$
described in Chapter~I of the book~\cite{GS1}. Let $\Omega$ and $M$
be monotone convex nonnegative differentiable indefinitely increasing
functions defined on the positive real semi-axis and satisfying the
condition $\Omega(0)=M(0)=0$. The space $W^\Omega_M(\Rset^k)$ is the
union (inductive limit) with respect to $A,B>0$ of the Banach spaces
consisting of entire analytic functions on $\Cset^k$ with the finite
norm
\[
\sup_{z=x+iy\in\Cset^k} |f(z)|\exp[M(Ax)-\Omega(By)].
\]
If $\Omega$ and $M$ grow faster than any linear function, then the
Fourier transformation isomorphically maps the space
$W^\Omega_M(\Rset^k)$ onto the space $W^{M_*}_{\Omega_*}(\Rset^k)$,
where
\[
M_*(s)=\sup_{t\geq 0} (st- M(t)),\quad
\Omega_*(s)=\sup_{t\geq 0} (st- \Omega(t))
\]
are the dual functions of $\Omega$ and $M$ in the sense of Young. For
$0 < \alpha\leq 1$ and $0\leq\beta<1$, the space
$S^\beta_\alpha(\Rset^k)$ coincides with the space
$W^\Omega_M(\Rset^k)$ with $\Omega(s)=s^{1/(1-\beta)}$ and
$M(s)=s^{1/\alpha}$. In particular,
$S^0_1(\Rset^k)=W^\Omega_M(\Rset^k)$, where $\Omega(s)=M(s)=s$. By
analogy with Definitions~\ref{dxx1} and~\ref{dxx2}, one can make

\begin{defn}
Let $U$ be a cone in $\Rset^k$ and $\Omega$ and $M$ be functions with
the properties specified above. The Banach space
$W^{\Omega,B}_{M,A}(U)$ consists of entire analytic functions on
$\Cset^k$ with the finite norm
\[
\sup_{z=x+iy\in \Cset^k} |f(z)|
\exp(M(|x/A|)-\Omega(\delta_U(Bx))-\Omega(|By|)),
\]
where $\delta_U(x)=\inf_{x'\in U}|x-x'|$ is the distance from $x$ to
$U$. The space $W^{\Omega}_M(U)$ is defined
by the relation
\[
W^\Omega_M(U)=\bigcup_{A,B>0,\,\tilde U\supset U}
W^{\Omega,B}_{M,A}(\tilde U),
\]
where $\tilde U$ runs over all conic neighborhoods of $U$
and the union is endowed with the
inductive limit topology.
\end{defn}

Further, we can
introduce the following definition analogous to
Definition~\ref{dxx3}:

\begin{defn}
Let $K$ be a closed cone in $\Rset^k$. The space $\mathcal
U^\Omega_M(K)$ is defined to be the inductive limit
$\varinjlim\nolimits_{K'\in \mathcal P(K)} W^{\prime \Omega}_M(K)$,
where $\mathcal P(K)$ is the set of all nonempty proper closed cones
contained in $K$. A closed cone $K$ is said to be a carrier cone of
an element $u\in \mathcal U^\Omega_M(\Rset^k)$ if the latter belongs
to the image of the canonical mapping from $\mathcal U^\Omega_M(K)$
to~$\mathcal U^\Omega_M(\Rset^k)$.
\end{defn}

The results obtained in this paper suggest the following conjecture:

\begin{hypo}
Let the defining functions $\Omega$ and $M$ be such that $M(s)\leq
\Omega(as)$ for some $a>0$. Then the following statements hold:

\begin{enumerate}
\item[$(1)$] The space $\mathcal U^\Omega_M(\Rset^k)$ is nontrivial
regardless of the triviality or nontriviality of
$W_{M}^{\Omega}(\Rset^k)$.

\item[$(2)$] If $W_{M}^{\Omega}(\Rset^k)$ is nontrivial, then
$\mathcal U^\Omega_M(\Rset^k)$ is canonically isomorphic to the space
$W_{M}^{\prime \Omega}(\Rset^k)$.

\item[$(3)$] Theorems~$\ref{txx3}$,~$\ref{txx4}$, and~$\ref{txx5}$
are
valid for the spaces $\mathcal U^\Omega_M(K)$.

\item[$(4)$] One can canonically define the Fourier transformation
that isomorphically maps $\mathcal U^\Omega_M(\Rset^k)$ onto
$\mathcal U_{\Omega_*}^{M_*}(\Rset^k)$.
\end{enumerate}
\end{hypo}

Note that the Fourier transformation on $\mathcal
U^\Omega_M(\Rset^k)$ cannot be constructed as that of
ultrafunctionals because the elements of $\mathcal
U^\Omega_M(\Rset^k)$ grow faster than exponentially and their Laplace
transformation is not well defined.

\def\theequation{\Alph{section}.\arabic{equation}}

\appendix
\section{Proof of Theorem~\ref{t1}}
\label{app1}

This appendix is organized as follows. We first introduce some
additional notation concerning inductive systems, which will also be
used in proving Lemma~\ref{l5} in Appendix~\ref{app2}. Then we derive
several auxiliary results (Lemmas~\ref{l1}--\ref{l4}) and, finally,
prove Theorem~\ref{t1}.

Let $\mathcal X$ be an inductive system over a
partially ordered set $A$. For $I\subset A$, we denote by
$T^{\mathcal X}_I$ the set of triples $(x,\alpha,\alpha')$ such that
$\alpha,\alpha'\in I$, $\alpha\leq\alpha'$, and $x\in \mathcal
X(\alpha)$. If $(x,\alpha,\alpha')\in T^{\mathcal X}_A$, then we set
$\sigma^{\mathcal X}(x,\alpha,\alpha')= \iota^{\mathcal X}_\alpha
x-\iota^{\mathcal X}_{\alpha'}\rho^{\mathcal X}_{\alpha\alpha'}x$
(recall that $\iota^{\mathcal X}_\alpha$ is the canonical embedding
of $\mathcal X(\alpha)$ into $\oplus_{\alpha'\in A} \mathcal
X(\alpha')$). We denote by $N^{\mathcal X}_I$ the subspace of
$\oplus_{\alpha'\in A} \mathcal X(\alpha')$ spanned by all
$\sigma^{\mathcal X}(x,\alpha,\alpha')$ with $(x,\alpha,\alpha')\in
T^{\mathcal X}_I$. For $I\subset A$, we denote by $M^{\mathcal X}_I$
the subspace $\oplus_{\alpha\in I} \mathcal X(\alpha)$ of the space
$\oplus_{\alpha\in A} \mathcal X(\alpha)$. Obviously, the space
$\varinjlim \mathcal X^I$ is isomorphic to $M^{\mathcal
X}_I/N^{\mathcal X}_I$. We denote by $j^{\mathcal X}_I$ the canonical
surjection from $M^{\mathcal X}_I$ onto $\varinjlim \mathcal X^I$. If
$I\subset J\subset A$, then we have
\begin{equation}\label{1}
\tau^{\mathcal X}_{I,\,J}j^{\mathcal X}_I x=
j^{\mathcal X}_J x,\quad x\in M^{\mathcal X}_I.
\end{equation}

We say that a subset $I$ of a partially ordered set $A$ is hereditary
if the relations $\alpha\in I$ and $\alpha'\leq \alpha$ imply that
$\alpha'\in I$.

\begin{lem}
\label{l1}
Let $\mathcal X$ be a prelocalizable inductive system of vector
spaces over a distributive quasi-lattice $A$. If $I$ is a hereditary
subset of $A$, then $N^{\mathcal X}_A\cap M^{\mathcal X}_I =
N^{\mathcal X}_I$.
\end{lem}

\begin{proof}
The inclusion $N^{\mathcal X}_I\subset N^{\mathcal X}_A\cap
M^{\mathcal X}_I$ is obvious. To prove the converse inclusion, it
suffices to show that $N^{\mathcal X}_J\cap M^{\mathcal X}_I \subset
N^{\mathcal X}_I$ for every finite $\wedge$-closed $J\subset A$. For
$\alpha\in J$, we denote by $k(\alpha)$ the cardinality $|J_\alpha|$
of the set $J_\alpha = \{\alpha'\in J\,|\, \alpha'\geq\alpha\}$. It
is obvious that $\alpha=\inf J_\alpha$. Therefore, if
$\alpha,\alpha'\in J$, $\alpha\ne\alpha'$, and $k(\alpha')\geq
k(\alpha)$, then we have $J_{\alpha}\ne J_{\alpha'}$ and,
consequently, $k(\alpha\wedge\alpha')=|J_{\alpha\wedge\alpha'}|\geq
|J_\alpha\cup J_{\alpha'}|>|J_\alpha|=k(\alpha)$. For $n\in \Nset$,
set $C_n =\{\alpha\in J\,|\, k(\alpha)\geq n\}$. We have
$J=C_1\supset C_2\supset\ldots \supset C_{|J|}= \{\tilde \alpha\}$,
where $\tilde \alpha=\inf J$, and $C_n=\varnothing$ for $n>|J|$. We
shall say that an $x\in N^{\mathcal X}_J\cap M^{\mathcal X}_I$ admits
a decomposition of order $n$ if there are a family of vectors
$x_{\alpha\alpha'}\in \mathcal X(\alpha)$ indexed by the set
$\{(\alpha,\alpha') : \alpha\in C_n,\, \alpha'\in J\setminus
I,\,\alpha<\alpha'\,\}$ and an element $\tilde x\in N^{\mathcal X}_I$
such that\footnote{Here and below, we assume that the sum of a family
of vectors indexed by the empty set is equal to zero.}
\begin{equation}
\label{2}
x= \tilde x + \sum_{\alpha\in C_n,\, \alpha'\in
J\setminus I,\,\alpha<\alpha'}\sigma^{\mathcal X}(x_{\alpha\alpha'},
\alpha,\alpha').
\end{equation}
If $x$ has a decomposition of order $>|J|$, then $x\in N^{\mathcal
X}_I$. Therefore, the lemma will be proved as soon as we show that
every $x\in N^{\mathcal X}_J\cap M^{\mathcal X}_I$ admits a
decomposition of order $n$ for any $n\in \Nset$. Since $I$ is
hereditary, every $x\in N^{\mathcal X}_J\cap M^{\mathcal X}_I$ has a
decomposition of order $1$, and we have to show that $x$ has a
decomposition of order $n+1$ supposing it has a decomposition of the
form~(\ref{2}) of order $n$. To this end, it suffices to establish
that $\sigma^{\mathcal X}(x_{\alpha\alpha'}, \alpha,\alpha')$ has a
decomposition of order $n+1$ for every $\alpha\in C_n,\, \alpha'\in
J\setminus I$ such that $\alpha<\alpha'$ and $k(\alpha)=n$. Let
$\Lambda=\{\beta\in C_n\,|\,\beta<\alpha',\,\beta\ne\alpha\}$. Since
$\alpha'\notin I$, the $\alpha'$-component of $x$ is equal to zero
and by (\ref{2}) we have
\begin{equation}\label{3}
\rho^{\mathcal X}_{\alpha\alpha'}\,x_{\alpha\alpha'}+
\sum_{\beta\in\Lambda}
\rho^{\mathcal X}_{\beta\alpha'}\,x_{\beta\alpha'}=0.
\end{equation}
If $\Lambda=\varnothing$, then the injectivity of $\rho^{\mathcal
X}_{\alpha\alpha'}$ implies that $x_{\alpha\alpha'}=0$ and
$\sigma^{\mathcal X}(x_{\alpha\alpha'}, \alpha,\alpha')=0$.
Therefore, in this case, $\sigma^{\mathcal
X}(x_{\alpha\alpha'},\alpha,\alpha')$ admits decompositions of all
orders. Now let $\Lambda\ne \varnothing$ and $\tilde
\beta=\sup\Lambda$ (the element $\tilde \beta$ is well defined
because $\Lambda$ is a finite set whose elements do not exceed
$\alpha'$; note that $\tilde \beta$ does not necessarily belong to
$J$). Set $y=\sum_{\beta\in \Lambda} \rho^{\mathcal
X}_{\beta\tilde\beta}\,x_{\beta\alpha'}$. Then it follows from
(\ref{3}) that $\rho^{\mathcal
X}_{\alpha\alpha'}\,x_{\alpha\alpha'}+\rho^{\mathcal
X}_{\tilde\beta\alpha'}\,y=0$. Hence, by III there is a $z\in
\mathcal X(\tilde\beta\wedge\alpha)$ such that
$x_{\alpha\alpha'}=\rho^{\mathcal
X}_{\tilde\beta\wedge\alpha,\,\alpha}\,z$. Because the quasi-lattice
$A$ is distributive, we have $\tilde\beta\wedge\alpha=\sup_{\beta\in
\Lambda} \beta\wedge\alpha$ and by II, there is a family
$\{z_\beta\}_{\beta\in\Lambda}$ such that $z_\beta\in \mathcal
X(\beta\wedge\alpha)$ and $z=\sum_{\beta\in \Lambda} \rho^{\mathcal
X}_{\beta\wedge\alpha,\,\tilde\beta\wedge\alpha}\,z_{\beta}$. We thus
have $x_{\alpha\alpha'}=\sum_{\beta\in \Lambda}\rho^{\mathcal
X}_{\beta\wedge\alpha,\,\alpha}\,z_{\beta}$ and, consequently,
\begin{equation}\label{4}
\sigma^{\mathcal X}(x_{\alpha\alpha'},\alpha,\alpha')=
\sum_{\beta\in \Lambda}
[\sigma^{\mathcal X}(z_{\beta},\alpha\wedge\beta,\alpha')-
\sigma^{\mathcal X}(z_{\beta},\alpha\wedge\beta,\alpha)].
\end{equation}
If $\alpha\in I$, then we set $\tilde y=- \sigma^{\mathcal
X}(z_{\beta},\alpha\wedge\beta,\alpha)$ and $y_{\gamma\gamma'}=
\delta_{\gamma',\,\alpha'}
\sum_{\beta\in\Lambda,\,\alpha\wedge\beta=\gamma}z_\beta$, where
$\gamma,\gamma'\in J$ and $\delta_{\gamma',\,\alpha'}=1$ for
$\gamma'=\alpha'$ and $\delta_{\gamma',\,\alpha'}=0$ for
$\gamma'\ne\alpha'$. If $\alpha\notin I$, then we set $\tilde y=0$,
$y_{\gamma\gamma'}= \delta_{\gamma',\,\alpha'}
\sum_{\beta\in\Lambda,\,\alpha\wedge\beta=
\gamma}z_\beta-\delta_{\gamma',\,\alpha}
\sum_{\beta\in\Lambda,\,\alpha\wedge\beta=\gamma}z_\beta$. Since
$k(\alpha\wedge\beta)>k(\alpha)=n$ for $\beta\in \Lambda$, it follows
from~(\ref{4}) that
\[
\sigma^{\mathcal X}(x_{\alpha\alpha'},\alpha,\alpha')=
\tilde y + \sum_{\gamma\in C_{n+1},\, \gamma'\in
J\setminus I,\,\gamma<\gamma'}\sigma^{\mathcal X}(y_{\gamma\gamma'},
\gamma,\gamma'),
\]
i.e., $\sigma^{\mathcal X}(x_{\alpha\alpha'},\alpha,\alpha')$ admits
a decomposition of order $n+1$. The lemma is proved.
\end{proof}

\begin{cor}
\label{c4}
Let $A$ be a distributive quasi-lattice, $\mathcal X$ be a
prelocalizable inductive system over $A$, and $I\subset J\subset A$.
If $I$ is a hereditary subset of $A$, then the canonical mapping
$\tau^{\mathcal X}_{I,\,J}\colon \varinjlim \mathcal X^I\to
\varinjlim \mathcal X^J$ is injective.
\end{cor}

\begin{proof}
Let $x\in \varinjlim \mathcal X^I$ and $\tau^{\mathcal
X}_{I,\,J}x=0$. By the surjectivity of $j^{\mathcal X}_I$, there is
an $\tilde x\in M^{\mathcal X}_I$ such that $x=j^{\mathcal X}_I
\tilde x$. It follows from (\ref{1}) that $j^{\mathcal X}_J \tilde
x=0$, i.e., $\tilde x\in N^{\mathcal X}_J$. Therefore, $\tilde x\in
M^{\mathcal X}_I\cap N^{\mathcal X}_J$ and in view of Lemma~\ref{l1}
we conclude that $\tilde x\in N^{\mathcal X}_I$ and $x=j^{\mathcal
X}_I \tilde x=0$. The corollary is proved.
\end{proof}

\begin{lem}
\label{l2}
Let $\mathcal X$ be an inductive system over a partially ordered set
$A$ and $I_1$ and $I_2$ be hereditary subsets of $A$. Then for every
$x\in N^{\mathcal X}_{I_1\cup I_2}$, there are $x_{1,2}\in
N^{\mathcal X}_{I_{1,2}}$ such that $x=x_1+x_2$.
\end{lem}

\begin{proof}
Let $\Lambda$ be the set of all pairs $(\alpha,\alpha')$ such that
$\alpha,\alpha'\in I_1\cup I_2$ and $\alpha \leq \alpha'$. By
definition of $N^{\mathcal X}_{I_1\cup I_2}$ there is a family
$\{x_{\alpha\alpha'}\}_{(\alpha,\alpha')\in \Lambda}$ such that
$x_{\alpha\alpha'}\in \mathcal X(\alpha)$ and
$x=\sum_{(\alpha,\alpha')\in \Lambda}\sigma^{\mathcal
X}(x_{\alpha\alpha'},\alpha,\alpha')$. We have $x=x_1+x_2$, where
\[
x_1=\sum_{(\alpha,\alpha')\in \Lambda,\,\alpha'\in I_1}
\sigma^{\mathcal X}(x_{\alpha\alpha'},\alpha,\alpha'), \quad
x_2=\sum_{(\alpha,\alpha')\in \Lambda,\,\alpha'\in I_2\setminus I_1}
\sigma^{\mathcal X}(x_{\alpha\alpha'},\alpha,\alpha').
\]
Since $I_{1,2}$ are hereditary, we conclude that $x_{1,2}\in
N^{\mathcal X}_{I_{1,2}}$. The lemma is proved.
\end{proof}

\begin{lem}
\label{l3}
Let $A$ be a distributive quasi-lattice, $\mathcal X$ be a
prelocalizable inductive system over $A$. Let $J\subset A$, and
$I_1,I_2$ be hereditary subsets of $A$ contained in $J$. Suppose
$x_{1,2}\in \varinjlim \mathcal X^{I_{1,2}}$ are such that
$\tau^{\mathcal X}_{I_1,\,J}x_1=\tau^{\mathcal X}_{I_2,\,J}x_2$. Then
there is an $x\in \varinjlim \mathcal X^{I_1\cap I_2}$ such that
$x_{1}=\tau^{\mathcal X}_{I_1\cap I_2,\,I_1}x$ and
$x_{2}=\tau^{\mathcal X}_{I_1\cap I_2,\,I_2}x$.
\end{lem}

\begin{proof}
Let $\tilde x_{1,2}\in M^{\mathcal X}_{I_{1,2}}$ be such that
$x_{1,2}=j^{\mathcal X}_{I_{1,2}}\tilde x_{1,2}$. We have
\[
\tau^{\mathcal X}_{I_1\cup I_2,\,J}
\tau^{\mathcal X}_{I_1,\,I_1\cup I_2}x_1=
\tau^{\mathcal X}_{I_1,\,J}x_1=\tau^{\mathcal X}_{I_2,\,J}x_2=
\tau^{\mathcal X}_{I_1\cup I_2,\,J}
\tau^{\mathcal X}_{I_2,\,I_1\cup I_2}x_2.
\]
Since the sets $I_{1,2}$ are hereditary, the set $I_1\cup I_2$ is
also hereditary and by Corollary~\ref{c4}, the mapping
$\tau^{\mathcal X}_{I_1\cup I_2,\,J}$ is injective. Therefore,
$\tau^{\mathcal X}_{I_1,\,I_1\cup I_2}x_1= \tau^{\mathcal
X}_{I_2,\,I_1\cup I_2}x_2$ and using (\ref{1}), we obtain
$j^{\mathcal X}_{I_1\cup I_2}(\tilde x_1- \tilde x_2)=0$. This means
that $\tilde x_1-\tilde x_2\in N^{\mathcal X}_{I_1\cup I_2}$. By
Lemma~\ref{l2}, there are $y_{1,2}\in N^{\mathcal X}_{I_{1,2}}$ such
that $\tilde x_1-\tilde x_2=y_1+y_2$. Set $\tilde x = \tilde x_1
-y_1=\tilde x_2+y_2$. Then $\tilde x\in M^{\mathcal X}_{I_1}\cap
M^{\mathcal X}_{I_2}= M^{\mathcal X}_{I_1\cap I_2}$. Set $x =
j^{\mathcal X}_{I_1\cap I_2}\tilde x$. Then $x\in \varinjlim \mathcal
X^{I_1\cap I_2}$ and it follows from (\ref{1}) that
\begin{align*}
&&\tau^{\mathcal X}_{I_1\cap I_2,\,I_1}x=
\tau^{\mathcal X}_{I_1\cap I_2,\,I_1}
j^{\mathcal X}_{I_1\cap I_2}\tilde x=
j^{\mathcal X}_{I_1}(\tilde x_1-y_1)=x_1, \\
&&\tau^{\mathcal X}_{I_1\cap I_2,\,I_2}x=
\tau^{\mathcal X}_{I_1\cap I_2,\,I_2}
j^{\mathcal X}_{I_1\cap I_2}\tilde x=
j^{\mathcal X}_{I_2}(\tilde x_2+y_2)=x_2.
\end{align*}
The lemma is proved.
\end{proof}

\begin{lem}
\label{l4}
Let $A$ be a quasi-lattice, $B$ be a lattice, and $\lambda\colon A\to
B$ be an injective quasi-lattice morphism such that any element
$\beta\in B$ is representable in the form
$\beta=\lambda(\alpha_1)\vee\ldots\vee\lambda(\alpha_n)$, where
$\alpha_1,\ldots,\alpha_n\in A$. Then we have

$1)$ If $\alpha,\alpha'\in A$ and
$\lambda(\alpha')\leq\lambda(\alpha)$, then $\alpha'\leq\alpha$.

$2)$ If $\beta,\beta'\in B$, $\beta'\leq \beta$, and
$\beta=\lambda(\alpha)$ for an $\alpha\in A$, then there is a unique
$\alpha'\in A$ such that $\beta'=\lambda(\alpha')$.

$3)$ If $A'\subset A$ has an infimum in $A$, then $\lambda(A')$ has
an infimum in $B$, and $\lambda(\inf A')=\inf \lambda(A')$.
\end{lem}

\begin{proof}
1) We have
$\lambda(\alpha\wedge\alpha')=
\lambda(\alpha)\wedge\lambda(\alpha')=\lambda(\alpha')$.
In view of the injectivity of $\lambda$ it hence follows that
$\alpha\wedge\alpha'=\alpha'$. This means that $\alpha'\leq\alpha$.

2) Let $\alpha_1,\ldots,\alpha_n\in A$ be such that
$\beta'=\lambda(\alpha_1)\vee\ldots\vee\lambda(\alpha_n)$. Since
$\lambda(\alpha_j)\leq \beta$, in view of~1) we have
$\alpha_j\leq\alpha$ for any $j=1,\ldots,n$. Therefore, the element
$\alpha'=\alpha_1\vee\ldots\vee\alpha_n$ is well defined and
satisfies the relation
$\lambda(\alpha')=
\lambda(\alpha_1)\vee\ldots\vee\lambda(\alpha_n)=\beta'$.
The uniqueness of $\alpha'$ follows from the injectivity of
$\lambda$.

3) Obviously, $\lambda(\inf A')\leq \beta'$ for any
$\beta'\in\lambda(A')$. Let $\beta\in B$ be such that
$\beta\leq\beta'$ for all $\beta'\in \lambda(A')$. Then by~2), there
is an $\alpha\in A$ such that $\beta=\lambda(\alpha)$, and in view
of~1) we have $\alpha\leq \alpha'$ for every $\alpha'\in A'$. This
implies that $\alpha\leq\inf A'$ and $\beta\leq \lambda(\inf A')$ and
so $\lambda(\inf A')=\inf \lambda(A')$.

The lemma is proved.
\end{proof}

\begin{proof}[Proof of Theorem~$\ref{t1}$.]
Let $\mathcal Z= \lambda(\mathcal X)$. Note that $A_\beta$ is a
hereditary subset of $A$ for any $\beta\in B$. The fulfilment of the
conditions I and III for $\mathcal Z$ therefore follows from
Corollary~\ref{c4} and from Lemma~\ref{l3} respectively. Let
$\beta_{1,2}\in B$, $\beta=\beta_1\vee \beta_2$, and $x\in \mathcal
Z(\beta)$. Since $\mathcal Z(\beta)=\varinjlim \mathcal
X^{A_{\beta}}$, there are $\alpha_1,\ldots,\alpha_m\in A$ and
$x_1\in\mathcal X(\alpha_1),\ldots,x_m\in\mathcal X(\alpha_m)$ such
that $\lambda(\alpha_j)\leq \beta$ and
$x=\sum_{j=1}^m\rho_{\alpha_j}^{\beta}\,x_j$, where
$\rho_{\alpha_j}^{\beta}$ is the canonical mapping from $\mathcal
X(\alpha_j)$ into $\varinjlim \mathcal X^{A_{\beta}}$. Choose
$\gamma_1^1,\ldots,\gamma_1^s,\gamma_2^1,\ldots,\gamma_2^t\in A$ such
that
\[
\beta_1=\lambda(\gamma_1^1)\vee\ldots\vee\lambda(\gamma_1^s),\quad
\beta_2=\lambda(\gamma_2^1)\vee\ldots\vee\lambda(\gamma_2^t).
\]
The distributivity of $B$ implies that
\[
\lambda(\alpha_j)=\lambda(\alpha_j)\wedge(\beta_1\vee\beta_2)=
\lambda((\alpha_j\wedge\gamma_1^1)\vee\ldots\vee (\alpha_j\wedge
\gamma_2^t)),\quad j=1,\ldots,m,
\]
and by the injectivity of $\lambda$, we have
$\alpha_j=(\alpha_j\wedge\gamma_1^1)\vee\ldots\vee
(\alpha_j\wedge\gamma_2^t)$. Since $\mathcal X$ satisfies the
condition~II, for any $j=1,\ldots,m$ there are $y^1_j\in \mathcal
X(\alpha_j\wedge\gamma_1^1),\ldots,y^s_j\in \mathcal
X(\alpha_j\wedge\gamma_1^s)$ and $z^1_j\in \mathcal
X(\alpha_j\wedge\gamma_2^1),\ldots,z^t_j\in \mathcal
X(\alpha_j\wedge\gamma_2^t)$ such that
\[
x_j=\sum_{l=1}^s \rho_{\alpha_j\wedge
\gamma_1^l,\,\alpha_j}^{\mathcal X}y^l_j+
\sum_{l=1}^t \rho_{\alpha_j\wedge
\gamma_2^l,\,\alpha_j}^{\mathcal X}z^l_j.
\]
Set $y=\sum_{j=1}^m\sum_{l=1}^s
\rho_{\alpha_j\wedge\gamma_1^l}^{\beta_1}y_j^l$,
$z=\sum_{j=1}^m\sum_{l=1}^t
\rho_{\alpha_j\wedge\gamma_2^l}^{\beta_2}z_j^l$. Then $y\in \mathcal
Z(\beta_1)$, $z\in \mathcal Z(\beta_2)$ and we have
\begin{multline}
\rho_{\beta_1,\,\beta}^{\mathcal Z}\,y+
\rho_{\beta_2,\,\beta}^{\mathcal Z}\,z=
\sum_{j=1}^m\sum_{l=1}^s
\rho_{\alpha_j\wedge\gamma_1^l}^{\beta}\,y_j^l +
\sum_{j=1}^m\sum_{l=1}^t
\rho_{\alpha_j\wedge\gamma_2^l}^{\beta}\,z_j^l=\\=
\sum_{j=1}^m \rho_{\alpha_j}^{\beta}\left[\sum_{l=1}^s
\rho_{\alpha_j\wedge\gamma_1^l,\,\alpha_j}^{\mathcal X}y^l_j+
\sum_{l=1}^t
\rho_{\alpha_j\wedge\gamma_2^l,\,\alpha_j}^{\mathcal X}z^l_j\right]=
\sum_{j=1}^m \rho_{\alpha_j}^{\beta}x_j=x.
\nonumber
\end{multline}
Thus, the inductive system $\mathcal X$ satisfies the condition~II
and, consequently, is prelocalizable.

We now suppose that the lattice $B$ is infinitely distributive and
$\mathcal X$ is a localizable inductive system and check that
$\mathcal Z$ satisfies the condition~III$'$. Let
$\{\beta_\omega\}_{\omega\in\Omega}$ be a nonempty family of elements
of $B$ bounded above by a $\beta\in B$, and let
$\{x_\omega\}_{\omega\in\Omega}$ be a family such that $x_\omega\in
\mathcal Z(\beta_\omega)$ and $y=\rho^{\mathcal
Z}_{\beta_\omega,\,\beta}\,x_\omega$ does not depend on $\omega$.

We first prove the statement for the case when
$\beta_{\omega_0}=\lambda(\alpha_0)$ for some $\omega_0\in\Omega$ and
$\alpha_0\in A$. For brevity, we write $\beta_0=\beta_{\omega_0}$ and
$x_0=x_{\omega_0}$. Set $\beta'_\omega= \beta_\omega\wedge \beta_0$.
Since $\beta'_\omega\leq\beta_0$, by Lemma~\ref{l4} there are
(uniquely defined) $\alpha_\omega\in A$ such that
$\alpha_\omega\leq\alpha_0$ and
$\beta'_\omega=\lambda(\alpha_\omega)$. Because $\mathcal Z$
satisfies the condition~III, there are $x'_\omega\in \mathcal
Z(\beta'_\omega)$ such that $\rho^{\mathcal
Z}_{\beta'_\omega,\,\beta_0}x'_\omega=x_0$ and $\rho^{\mathcal
Z}_{\beta'_\omega,\,\beta_\omega}x'_\omega=x_\omega$ for every
$\omega\in\Omega$. The canonical mapping
$\rho_{\alpha}^{\lambda(\alpha)}$ from $\mathcal X(\alpha)$ into
$\mathcal Z(\lambda(\alpha))=\varinjlim \mathcal
X^{A_{\lambda(\alpha)}}$ is isomorphic for any $\alpha\in A$ because
$\lambda(\alpha)$ is the biggest element of the set
$A_{\lambda(\alpha)}$. Therefore, for any $\omega\in \Omega$ there
exists a unique $\tilde x_\omega\in\mathcal X(\alpha_\omega)$ such
that $x'_\omega=\rho_{\alpha_\omega}^{\beta'_\omega}\tilde x_\omega$.
We have $\rho^{\mathcal
Z}_{\lambda(\alpha'),\,\lambda(\alpha)}
\rho^{\lambda(\alpha')}_{\alpha'}=
\rho^{\lambda(\alpha)}_\alpha\rho^{\mathcal X}_{\alpha'\alpha}$ for
any $\alpha',\alpha\in A$ such that $\alpha'\leq \alpha$. Hence
$\rho_{\alpha_0}^{\beta_0}\rho^{\mathcal
X}_{\alpha_\omega,\,\alpha_0}\tilde x_\omega= \rho^{\mathcal
Z}_{\beta'_\omega,\beta_0}x'_\omega=x_0$ and, consequently,
$\rho^{\mathcal X}_{\alpha_\omega,\,\alpha_0}\tilde
x_\omega=\left(\rho_{\alpha_0}^{\beta_0}\right)^{-1}x_0$ does not
depend on $\omega$. Let
$\tilde\alpha=\inf_{\omega\in\Omega}\alpha_\omega$ and
$\tilde\beta=\lambda(\tilde\alpha)$. By Lemma~\ref{l4}, we have
$\tilde\beta=\inf_{\omega\in\Omega}\beta'_\omega=
\inf_{\omega\in\Omega}\beta_\omega$.
In view of the localizability of $\mathcal X$ there is an $\tilde
x\in \mathcal X(\tilde \alpha)$ such that $\tilde
x_\omega=\rho^{\mathcal X}_{\tilde\alpha,\,\alpha_\omega}\tilde x$
for all $\omega\in\Omega$. Set
$x=\rho_{\tilde\alpha}^{\tilde\beta}\tilde x$. Then $x\in\mathcal
Z(\tilde\beta)$ and we have
\[
\rho^{\mathcal Z}_{\tilde\beta,\,\beta_\omega}x=
\rho^{\mathcal Z}_{\beta'_\omega,\,\beta_\omega}
\rho^{\mathcal Z}_{\tilde\beta,\,\beta'_\omega}
\rho_{\tilde\alpha}^{\tilde\beta}\tilde x=
\rho^{\mathcal Z}_{\beta'_\omega,\,\beta_\omega}
\rho_{\alpha_\omega}^{\beta'_\omega}
\rho^{\mathcal X}_{\tilde\alpha,\,\alpha_\omega}\tilde x=
\rho^{\mathcal Z}_{\beta'_\omega,\,\beta_\omega}x'_\omega=x_\omega.
\]

We now consider the general case. Let $\tilde
\beta=\inf_{\omega\in\Omega}\beta_\omega$ and $J$ be a finite
$\wedge$-closed subset of $B$ such that $J\subset \lambda(A)$ and
$\beta=\sup_{\alpha\in J}\lambda(\alpha)$. As in the proof of
Lemma~\ref{1}, we denote by $J_\gamma$ ($\gamma\in J$) the set
$\{\gamma'\in J\,|\, \gamma'\geq\gamma\}$. For $n\in \Nset$, set $C_n
=\{\gamma\in J\,|\, |J_\gamma|\geq n\}$. We have $J=C_1\supset
C_2\supset\ldots \supset C_{|J|}= \{\tilde \gamma\}$, where $\tilde
\gamma=\inf J$, and $C_n=\varnothing$ for $n>|J|$. It suffices to
show that for any $n\in \Nset$, there is a family
$\{y_\gamma\}_{\gamma\in C_n}$ such that $y_\gamma\in \mathcal
Z(\gamma)$ and
\begin{equation}\label{5}
y=\rho^{\mathcal Z}_{\tilde\beta,\,\beta}\,\tilde y+
\sum_{\gamma\in C_n}\rho^{\mathcal Z}_{\gamma,\,\beta}\,y_{\gamma},
\end{equation}
where $\tilde y\in \mathcal Z(\tilde\beta)$. We prove this statement
by induction on $n$. For $n=1$, the existence of a decomposition of
the form~(\ref{5}) follows from the condition~II. Therefore, it
suffices to show that if~(\ref{5}) holds for some $n\in \Nset$, then
for any $\gamma\in C_n$ there is a family
$\{y'_{\gamma'}\}_{\gamma'\in C_{n+1}}$ such that $y_{\gamma'}\in
\mathcal Z(\gamma')$ and
\begin{equation}\label{6}
y_\gamma=
\rho^{\mathcal Z}_{\tilde\beta\wedge\gamma,\,\gamma}\,
\tilde y_\gamma+
\sum_{\gamma'\in C_{n+1}}
\rho^{\mathcal Z}_{\gamma',\,\gamma}\,y'_{\gamma'},
\end{equation}
where $\tilde y_\gamma\in \mathcal Z(\tilde\beta\wedge\gamma)$. Let
$\Omega'$ be the disjoint union of $\Omega$ and a one-element set
$\{\varkappa\}$ ($\varkappa\notin\Omega$). Set
$\beta'_\varkappa=\gamma$ and $\beta'_\omega=\beta_\omega\vee\sup
(C_{n}\setminus\{\gamma\})$ for $\omega\in\Omega$ (if
$C_n\setminus\{\gamma\}=\varnothing$, then we assume
$\beta'_\omega=\beta_\omega$). For every $\omega\in \Omega'$, we
define an element $x'_\omega\in \mathcal Z(\beta'_\omega)$ setting
$x'_\varkappa=y_\gamma$ and
\[
x'_\omega=\rho^{\mathcal Z}_{\beta_\omega,\,\beta'_\omega}\,x_\omega-
\rho^{\mathcal Z}_{\tilde\beta,\,\beta'_\omega}\,\tilde y -
\sum_{\gamma'\in C_n\setminus\{\gamma\}}
\rho^{\mathcal Z}_{\gamma',\,\beta'_\omega}\,y_{\gamma'}
\]
for $\omega\in\Omega$. It follows from~(\ref{5}) that the element
$\rho^{\mathcal Z}_{\beta'_\omega,\,\beta}\,x'_\omega$ does not
depend on $\omega\in\Omega'$. Let
$\beta'=\inf_{\omega\in\Omega'}\beta'_\omega$. Since
$\beta'_\varkappa\in\lambda(A)$, we can apply the result of the
preceding paragraph and find an $x'\in \mathcal Z(\beta')$ such that
$y_\gamma=x'_\varkappa=\rho^{\mathcal Z}_{\beta',\,\gamma}\,x'$.
Because the lattice $B$ is infinitely distributive, we have
\[
\beta'=(\tilde\beta\wedge\gamma)\vee
\sup_{\gamma'\in C_{n}\setminus\{\gamma\}} (\gamma'\wedge\gamma)
\quad (\beta'=\tilde\beta\wedge\gamma \mbox{ for }
C_n\setminus\{\gamma\}=\varnothing).
\]
Hence, by~II, there are $\tilde y_\gamma\in \mathcal
Z(\tilde\beta\wedge\gamma)$ and a family $\{z_{\gamma'}\}_{\gamma'\in
C_n\setminus\{\gamma\}}$ such that $z_{\gamma'}\in \mathcal
Z(\gamma'\wedge\gamma)$ and $y_\gamma=\rho^{\mathcal
Z}_{\tilde\beta\wedge\gamma,\,\gamma}\,\tilde y_\gamma +
\sum_{\gamma'\in C_n\setminus\{\gamma\}}\rho^{\mathcal
Z}_{\gamma'\wedge\gamma,\,\gamma}\,y'_\gamma$. Because
$\gamma'\wedge\gamma\in C_{n+1}$ for $\gamma',\gamma\in C_n$ and
$\gamma'\ne \gamma$, we can rewrite the last decomposition in the
form~(\ref{6}). Theorem~\ref{t1} is proved.
\end{proof}

\section{Proof of Lemma~\ref{l5}}
\label{app2}

In what follows, we use the notation introduced in the beginning
of Appendix~\ref{app1}.

Let $l$ be the linear mapping from $\oplus_{\tau\in T} \mathcal
X(\lambda(\tau))$ to $\oplus_{\alpha\in A}\mathcal X(\alpha)$ such
that $l j_\tau=\iota^{\mathcal X}_{\lambda(\tau)}$ for any $\tau\in
T$. The operator $l$ carries the vector~(\ref{11}) to the element
\begin{equation}\label{11aa}
\iota^{\mathcal X}_{\alpha}
\rho^{\mathcal X}_{\alpha\wedge\alpha',\,\alpha} x-
\iota^{\mathcal X}_{\alpha'}
\rho^{\mathcal X}_{\alpha\wedge\alpha',\,\alpha'} x=
\sigma^{\mathcal X}(x, \alpha\wedge\alpha',\alpha')-
\sigma^{\mathcal X}(x, \alpha\wedge\alpha',\alpha),
\end{equation}
where $\alpha=\lambda(\tau)$ and $\alpha'=\lambda(\tau')$. This
implies that $l(\EuScript N)\subset N^{\mathcal X}_A$ and hence
$\EuScript N\subset\mathrm{Ker}\,j^{\mathcal X}l$. The mapping
$j^{\mathcal X}l$ therefore uniquely determines a mapping $m\colon
\oplus_{\tau\in T} \mathcal X(\lambda(\tau))/\EuScript N\to
\varinjlim \mathcal X$. To prove the lemma, we have to show that $m$
is an isomorphism. To this end, it suffices to establish the opposite
inclusion
\begin{equation}\label{11a}
\EuScript N\supset\mathrm{Ker}\,j^{\mathcal X}l.
\end{equation}
Set $I=\lambda(T)$.
Let a mapping $\lambda' \colon I\to T$ be such that
$\lambda(\lambda'(\alpha))=\alpha$ for any $\alpha\in I$ and let the
mapping $l'\colon M_I^{\mathcal X}\to \oplus_{\tau\in T} \mathcal
X(\lambda(\tau))$ be defined by the relations $l'\iota^{\mathcal
X}_{\alpha}x=j_{\lambda'(\alpha)}x$ for any $\alpha\in I$ and
$x\in\mathcal X(\alpha)$. Clearly, $\mathrm{Im}\,l'$ coincides with
the subspace $E$ of $\oplus_{\tau\in T} \mathcal X(\lambda(\tau))$
spanned by all elements $j_\tau x$ with $\tau\in\lambda'(I)$ and
$x\in \mathcal X(\lambda(\tau))$. Moreover, we have $l'lx=x$ for any
$x\in E$. Let us show that every $x\in \oplus_{\tau\in T} \mathcal
X(\lambda(\tau))$ can be decomposed as $x=n+x'$, where $n\in\EuScript
N$ and $x'\in E$. It suffices to consider the case $x=j_\tau y$,
where $\tau\in T$ and $y\in \mathcal X(\lambda(\tau))$. Let
$\tau'=\lambda'(\lambda(\tau))$. Then we have
$\lambda(\tau)=\lambda(\tau')$ and, consequently, the element
$n=j_\tau y-j_{\tau'} y$ belongs to $\EuScript N$. Setting
$x'=x-n=j_{\tau'} y$, we obtain the desired decomposition because
$\tau'\in \lambda'(I)$.

Let $\tilde {\EuScript N}$ be the subspace of $M_I^{\mathcal X}$
spanned by all vectors of the form $\iota^{\mathcal X}_{\alpha}
\rho^{\mathcal X}_{\alpha\wedge\alpha',\,\alpha} x- \iota^{\mathcal
X}_{\alpha'} \rho^{\mathcal X}_{\alpha\wedge\alpha',\,\alpha'} x$
with $\alpha,\alpha'\in I$ and $x\in \mathcal
X(\alpha\wedge\alpha')$. We obviously have
\begin{equation}\label{11bb}
l(\EuScript N)=\tilde{\EuScript N},\quad l'(\tilde{\EuScript N})
\subset \EuScript N.
\end{equation}
The inclusion~(\ref{11a}) can be easily derived from the equality
\begin{equation}\label{11b}
N^{\mathcal X}_A\cap M_I^{\mathcal X}=\tilde {\EuScript N}
\end{equation}
which will be proved a little bit later. Indeed, let $x\in
\mathrm{Ker}\,j^{\mathcal X}l$. Then we have $lx\in
\mathrm{Ker}\,j^{\mathcal X}= N^{\mathcal X}_A$ and in view of the
obvious inclusion $\mathrm{Im}\,l\subset M^{\mathcal X}_I$ it follows
from~(\ref{11b}) that $lx\in\tilde{\EuScript N}$. According to the
above we can write $x=n+x'$, where $n\in \EuScript N$ and $x'\in E$.
By~(\ref{11bb}), we have $ln\in\tilde {\EuScript N}$ and, therefore,
$lx'\in \tilde{\EuScript N}$. Since $x'\in E$, we have $x'=l'lx'$ and
it follows from~(\ref{11bb}) that $x'\in\EuScript N$. Thus,
$x\in\EuScript N$ and the implication
(\ref{11b})~$\Rightarrow$~(\ref{11a}) is proved.

It remains to prove~(\ref{11b}). The inclusion $\tilde {\EuScript
N}\subset N^{\mathcal X}_A\cap M_I^{\mathcal X}$ obviously follows
from~(\ref{11aa}) and we have to verify that $x\in\tilde{\EuScript
N}$ supposing $x\in N^{\mathcal X}_A\cap M^{\mathcal X}_I$. Let
$\alpha',\alpha\in A$ be such that $\alpha'\leq \alpha$ and let $y\in
\mathcal X(\alpha')$. Since the set $I$ is cofinal in $A$, there is a
$\beta\in I$ such that $\beta\geq \alpha$, and we have
$\sigma^{\mathcal X}(y, \alpha',\alpha)=\sigma^{\mathcal X}(y,
\alpha',\beta)- \sigma^{\mathcal X}(\rho^{\mathcal X}_{\alpha'\alpha}
y, \alpha,\beta)$. Therefore, when writing sums of the elements of
the form $\sigma^{\mathcal X}(y, \alpha',\alpha)$, we can always
assume that $\alpha\in I$. In particular, since $x\in N^{\mathcal
X}_A$, we can write
\begin{multline}
x = \sum_{(\alpha',\alpha)\in A\times I,\,\alpha'\leq\alpha}
\sigma^{\mathcal X}(x_{\alpha'\alpha}, \alpha',\alpha)=\\
\sum_{\alpha'\in A\setminus I}
\sum_{\alpha\in C(\alpha')}
\sigma^{\mathcal X}(x_{\alpha'\alpha}, \alpha',\alpha)+
\sum_{(\alpha',\alpha)\in I\times I,\,\alpha'\leq\alpha}
\sigma^{\mathcal X}(x_{\alpha'\alpha}, \alpha',\alpha),
\label{11c}
\end{multline}
where $C(\alpha')=\{\,\alpha\in I\,|\, \alpha'\leq \alpha\,\}$ and
the family $\{x_{\alpha'\alpha}\}_{(\alpha',\alpha)\in A\times I}$
contains only finite number of nonzero elements. It is obvious that
the second sum in the right-hand side belongs to $\tilde {\EuScript
N}$. Therefore, it suffices to show that $y_{\alpha'}=\sum_{\alpha\in
C(\alpha')} \sigma^{\mathcal X}(x_{\alpha'\alpha}, \alpha',\alpha)$
belongs to $\tilde {\EuScript N}$ for any given $\alpha'\in
A\setminus I$. Since the $\alpha'$-component of $x$ is equal to zero,
the equality~(\ref{11c}) implies that $\sum_{\alpha\in C(\alpha')}
x_{\alpha'\alpha}=0$. Fixing an $\tilde \alpha\in C(\alpha')$, we
therefore obtain
\[
y_{\alpha'}=\sum_{\alpha\in C(\alpha')\setminus\{\tilde\alpha\}}
(\sigma^{\mathcal X}(x_{\alpha'\alpha}, \alpha',\alpha)-
\sigma^{\mathcal X}(x_{\alpha'\alpha}, \alpha',\tilde\alpha)).
\]
Using~(\ref{11aa}) and the relation
\[
\sigma^{\mathcal X}(x_{\alpha'\alpha}\,, \alpha',\alpha)-
\sigma^{\mathcal X}(x_{\alpha'\alpha}\,, \alpha',\tilde\alpha)=
\sigma^{\mathcal X}(z, \alpha\wedge\tilde\alpha,\alpha)-
\sigma^{\mathcal X}(z, \alpha\wedge\tilde\alpha,\tilde\alpha),
\]
where $z= \rho^{\mathcal X}_{\alpha',\,\alpha\wedge\tilde\alpha}
x_{\alpha'\alpha}$, we conclude that
$y_{\alpha'}\in \tilde{\EuScript N}$. The lemma is
proved.


\begin{thebibliography}{30}

\bibitem{GS}  I.~M.~Gelfand, G.~E.~Shilov, Generalized functions,
Vol. 2: Spaces of fundamental and generalized functions, Academic
Press, New York--London, 1968.

\bibitem{GS1} I.~M.~Gelfand, G.~E.~Shilov, Generalized functions,
Vol. 3: Theory of differential equations, Academic Press, New
York--London, 1967.

\bibitem{Graetzer} G.~Gr\"atzer, General Lattice Theory,
Akademie-Verlag, Berlin, 1978.

\bibitem{Hoermander} L.~H\"{o}rmander, The Analysis of Linear Partial
Differential Operators, Vol. I, Springer-Verlag,
Berlin--Heidelberg--New York--Tokyo, 1983.

\bibitem{Kawai} T.~Kawai, On the theory of Fourier hyperfunctions and
its applications to partial differential equations with constant
coefficients, J. Fac. Sci. Univ. Tokio (Sec.~1A Math.) 17 (1970)
467--517.

\bibitem{Komatsu} H.~Komatsu, Projective and injective limits of
weakly compact sequences of locally convex spaces, J.~Math. Soc.
Japan 19 (1967) 366--383.

\bibitem{Komatsu1} H.~Komatsu, Ultradistributions, I. Structure
theorems and a characterization, J.~Fac. Sci. Univ. Tokio (Sec. 1A
Math.) 20 (1973) 25--105.

\bibitem{Komatsu2} H.~Komatsu, An introduction to the theory of
hyperfunctions, in: Hyperfunctions and Pseudo-Differential equations,
Lecture Notes in Math., Vol.~287, Springer-Verlag,
Berlin--Heidelberg--New York, 1973, pp.~3--40.

\bibitem{Morimoto} M.~Morimoto, An Introduction to Sato's
Hyperfunctions, American Mathematical Society, Providence, Rhode
Island, 1993.

\bibitem{R}  C.~Roumieu, Sur quelques extensions de la notion de
distribution, Ann. Sci. \'Ecole Norm. Sup. (Ser.~3) 77 (1960)
41--121.

\bibitem{Sato} M.~Sato, Theory of hyperfunctions, S\^ugaku 10 (1958)
1--27 (Japanese).

\bibitem{Smirnov} A.~G.~Smirnov, Towards Euclidean theory of infrared
singular quantum fields, J.~Math. Phys. 44 (2003) 2058--2076.

\bibitem{Sol1} M.~A.~Soloviev, Towards a generalized distribution
formalism for gauge quantum fields, Lett. Math. Phys. 33 (1995)
49--59.

\bibitem{Sol2} M.~A.~Soloviev, An extension of distribution theory
and of the Paley--Wiener--Schwartz theorem related to quantum gauge
theory, Commun. Math. Phys. 184 (1997) 579--596.

\bibitem{Sol3} M.~A.~Soloviev, PCT, spin and statistics and analytic
wave front set, Theor. Math. Phys. 121 (1999) 1377--1396.

\end{thebibliography}
\end{document}